\documentclass[review, 3p, compress]{elsarticle}
\usepackage{lineno,hyperref}
\usepackage{amsfonts}
\usepackage{amssymb}
\usepackage{amsmath}
\usepackage{cases}
\usepackage{graphicx}
\usepackage{rotating}
\usepackage{float}
\usepackage{booktabs}
\usepackage{multirow,multicol}
\usepackage{epstopdf}
\usepackage{subfigure}
\usepackage{color}
\usepackage{setspace}
\usepackage{bm}
\usepackage{doi}
\usepackage{algorithm,algorithmic}
\usepackage{natbib}
\usepackage{subeqnarray}

\makeatletter
\@addtoreset{equation}{section}
\makeatother

\newtheorem{theorem}{Theorem}[section]
\newtheorem{lemma}{Lemma}[section]
\newtheorem{remark}{Remark}
\modulolinenumbers[5]

\journal{JSC}

\begin{document}

\begin{frontmatter}

\title{A preconditioning technique for an all-at-once system 
	from Volterra subdiffusion equations with graded time steps}

\author[address1]{Yong-Liang Zhao}%\corref{correspondingauthor}
\ead{ylzhaofde@sina.com}

\author[address2]{Xian-Ming Gu\corref{correspondingauthor}}
%\cortext[correspondingauthor]{Corresponding author}
\ead{guxianming@live.cn}

\author[address3]{Alexander Ostermann\corref{correspondingauthor}}
\ead{alexander.ostermann@uibk.ac.at}
\cortext[correspondingauthor]{Corresponding authors}

%\author[address1]{Pei-Yong Zhu}
%\ead{zpy6940@uestc.edu.cn}

%\author[address1]{Xi-Le Zhao}
%\ead{xlzhao122003@163.com}
%
%\author[address1]{Huan-Yan Jian}
%\ead{uestc_hyjian@sina.com}

\address[address1]{School of Mathematical Sciences, \\
University of Electronic Science and Technology of China, \\
Chengdu, Sichuan 611731, P.R. China}
\address[address2]{School of Economic Mathematics/Institute of Mathematics,\\
	Southwestern University of Finance and Economics,\\
	Chengdu, Sichuan 611130, P.R. China}
\address[address3]{Department of Mathematics, University of Innsbruck, \\
Technikerstra{\ss}e 13, Innsbruck 6020, Austria}
%\address[address3]{School of Science, Lanzhou University of Technology, \\
%Lanzhou, Gansu 730050, P.R. China}

\begin{abstract}
Volterra subdiffusion problems with weakly singular kernel describe the dynamics of subdiffusion processes well.
The graded $L1$ scheme is often chosen to discretize such problems 
since it can handle the singularity of the solution near $t = 0$.
In this paper, we propose a modification. We first split the time interval $[0, T]$ into $[0, T_0]$ and $[T_0, T]$, 
where $T_0$ ($0 < T_0 < T$) is reasonably small.
Then, the graded $L1$ scheme is applied in $[0, T_0]$, while the uniform one is used in $[T_0, T]$.
Our all-at-once system is derived based on this strategy. In order to solve the arising system efficiently,
we split it into two subproblems and design two preconditioners.
Some properties of these two preconditioners are also investigated.
Moreover, we extend our method to solve semilinear subdiffusion problems.
Numerical results are reported to show the efficiency of our method.
\end{abstract}


\end{frontmatter}

%\linenumbers

\section{Introduction}
\label{sec1}

Anomalous diffusion phenomena are common in various complex systems such as amorphous semiconductors~\cite{sokolov2002fractional},
filled polymers~\cite{metzler1995relaxation}, porous systems~\cite{he1998approximate} and turbulent plasma~\cite{del2004fractional}.
Generally, they do not have the Markovian stochastic property and cannot be simulated by the classical diffusion equations arising from Fick's law.
One of the widely used approaches for modelling anomalous diffusion are fractional diffusion equations \cite{metzler2000boundary,gorenflo2002time,podlubny2009matrix,pagnini2016stochastic}.
Numerous studies about fractional partial differential equations can be found in
\cite{podlubny1998fractional,lin2019crank,lei2018fast,shen2020h2n2,liao2019discrete,cao2019blow,gu2020parallel,zhao2020preconditioning,li2018fast}
and the references therein.

In this work, we are interested in solving the subdiffusion equation:
\begin{equation}
\label{eq1.1}
\begin{cases}
\int_{0}^{t} \xi_{1 - \beta}(t - s) \partial_{s} u(x,y,s) ds = \kappa \Delta u(x,y,t) + f(x,y,t), 
& (x,y,t) \in \Omega \times (0,T], \\
u(x,y,t) = 0, & (x,y) \in \partial \Omega,~0 < t \leq T, \\
u(x,y,0) = u_0(x,y), & (x,y) \in \Omega,
\end{cases}
\end{equation}
where $0 < \beta < 1$, the diffusion coefficient $\kappa > 0$, 
$\Omega = [x_L, x_R] \times [y_L,y_R] \subset \mathbb{R}^{2}$ and
the weakly singular kernel $\xi_\gamma(t) = t^{\gamma - 1}/\Gamma(\gamma)~(t > 0)$.

Eq.~\eqref{eq1.1} describes the dynamics of subdiffusion processes well,
in which the mean square variance grows at a rate slower than in a Gaussian process \cite{bouchaud1990anomalous}.
For the analytical solution of Volterra problems with weakly singular kernel (so-called time-fractional problems),
the Fourier transform method, the Laplace transform method and the Mellin transform method can be used \cite{podlubny1998fractional}.
However, in real applications, these methods are inappropriate for most time-fractional problems
because of the nonlocality and complexity of the fractional derivatives.
Thus, the development of efficient and reliable numerical techniques for solving Eq.~\eqref{eq1.1} attracts many researchers.

Until now, numerous articles have been published for solving time-fractional problems numerically in an efficient way
\cite{lin2007finite,li2009space,gao2014new,zhang2014finite,alikhanov2015new,jin2016analysis,zeng2015numerical,hu2020using}.
In these studies, the $L1$ or $L1$-type approximation is the most considered method to approach the Volterra operator in Eq.~\eqref{eq1.1}.
The convergence rate of this approximation on a uniform mesh is $2 - \beta$ under the assumption that $u$ is smooth on the closure of the domain \cite{lin2007finite}.
This smoothness assumption is unrealistic and ignores the weak singularity near the initial time $t = 0$.
This singularity has a great influence on the convergence rate.
In order to compensate for the singular behaviour of the exact solution at $t = 0$,
Mustapha et al.~\cite{mustapha2012finite,mustapha2011implicit}
applied the $L1$ formula with non-uniform time step to solve a class of subdiffusion equations.
Zhang et al.~\cite{zhang2014finite} investigated the non-uniform $L1$ approximation and applied it to solve fractional diffusion equations.
Stynes et al.~\cite{stynes2017error} proposed a finite difference scheme on a graded mesh in time to numerically solve time-fractional diffusion equations.
Liao et al.~\cite{liao2018sharp} analyzed the non-uniform $L1$ approximation for solving reaction-subdiffusion equations.
Other techniques for improving the poor accuracy of numerical approaches for the Volterra operator in Eq.~\eqref{eq1.1} can be found in \cite{jin2016analysis,zeng2015numerical,lubich1996nonsmooth,zeng2017second,yan2018analysis,
jin2019numerical,wang2020higher}.

The numerical schemes mentioned above are time-stepping schemes.
This means that the numerical solutions are obtained step-by-step.
Another class of methods consists in computing the numerical solutions in a parallel-in-time (PinT) pattern.
This includes the Laplacian inversion technique \cite{kwon2003parallel,mclean2010maximum} 
and the parareal algorithm \cite{li2013parallel,wu2018parareal,fu2019preconditioned}.
Recently, all-at-once systems arising from fractional partial differential equations have been studied by many researchers \cite{gu2020parallel,zhao2020preconditioning,ke2015fast,huang2017fast,lu2015fast,lu2018approximate,
Bertaccini2018limited,BERTACCINI201992}. 
Ke et al.~\cite{ke2015fast} developed a fast solver based on the divide-and-conquer (DC) method for the all-at-once system arising from time-fractional partial differential equations.
Lu et al.~\cite{lu2015fast} proposed an approximate inversion (AI) method to solve the block lower triangular Toeplitz linear system with tridiagonal blocks from fractional sub-diffusion equations.
Later, according to the short-memory principle \cite{Bertaccini2017mixed},
Bertaccini and Durastante \cite{Bertaccini2018limited} proposed a limited memory block preconditioner 
for fast solving their linear systems from space-fractional equations. 
In \cite{BERTACCINI201992}, the authors studied an all-at-once system arising from high-dimensional space-fractional equations.
In order to solve their system efficiently, they rewrote it into a tensor form 
and designed a tensor structured preconditioner.
However, in these articles, fast algorithms are designed based on uniform meshes.
In this paper, we try to get the numerical solution of Eq.~\eqref{eq1.1} globally in time by solving the corresponding all-at-once system with variable time steps.
In~\cite{mustapha2012finite,mustapha2011implicit,stynes2017error,liao2018sharp}, the authors used a time graded mesh when approximating the Caputo fractional derivative.
However, if the temporal regularity of the solution is small, the grid points of the time graded mesh become very dense near $t = 0$ and very sparse near $t = T$.
This will reduce the numerical resolution of the solution. It is also a bad choice to derive an all-at-once system based on such a non-uniform mesh,
because the storage requirement would be terrible and cannot be reduced.
The authors in \cite[Remark 8]{liao2018sharp} suggested that if $u$ is smooth away from $t = 0$,
the time interval $[0,T]$ can be split into $[0,T_0]$ and $[T_0,T]$, where $T_0$ ($0 < T_0 < T$) is reasonably small.
Then, the graded mesh is applied to $[0,T_0]$ and the uniform mesh is used in $[T_0,T]$.
According to their suggestion, we derive an all-at-once system in this paper based on such a strategy.
Further, a fast algorithm is designed to solve this system efficiently.

The rest of this paper is organized as follows. Section \ref{sec2} derives our all-at-once system of Eq.~\eqref{eq1.1}.
In Section \ref{sec3}, the system is split into two subproblems,
and two preconditioners are designed to efficiently solve these subproblems.
Moreover, several properties of these preconditioners are investigated.
In Section \ref{sec4}, we extend our algorithm to solve the semilinear problem of Eq.~\eqref{eq1.1}.
In Section \ref{sec5}, numerical results are reported. Concluding remarks are given in Section \ref{sec6}.

\section{The all-at-once system}
\label{sec2}

In this section, we first derive the time-stepping scheme for approximating Eq.~\eqref{eq1.1} 
by using the finite difference method.
Then, the all-at-once system is obtained based on this scheme.

\subsection{The time-stepping scheme}
\label{sec2.1}

For a given reasonably small $T_0~(0 < T_0 < T)$, we split $[0,T]$ into two parts $[0,T_0]$ and $[T_0,T]$.
In the first part $[0,T_0]$, we use the graded mesh $t_k = T_0 \left( \frac{k}{M_0} \right)^r$ for $k = 0, 1,\cdots, M_0$,
where $M_0$ is a positive integer and $r \geq 1$ is the grading parameter.
In the second part $[T_0,T]$, the uniform mesh $t_{M_0 + k} = T_0 + k \tilde{\tau}~(k = 1,2,\cdots,M - M_0)$ is used with $\tilde{\tau} = \frac{T - T_0}{M - M_0}$,
where $M > M_0$ is a positive integer.
Then, we get the mesh points $\left\{ t_k \right\}_{k = 0}^{M}$ and denote the time steps $\tau_k = t_k - t_{k - 1}~(k = 1, 2, \cdots, M)$.
Notice that for $M_0 + 1 \leq k \leq M$, we have $\tau_k = \tilde{\tau}$.
Let $h_x = \frac{x_R - x_L}{N_x}$ and $h_y = \frac{y_R - y_L}{N_y}$ 
be the grid spacing in $x$ and $y$ directions for given positive integers $N_x$ and $N_y$. Hence the space domain is discretized
by $\bar{\omega}_{h} = \left\{ (x_i, y_j) = (x_L + i h_x, y_L + j h_y) 
\mid 0 \leq i \leq N_x,~0 \leq j \leq N_y \right\}$.
According to \cite{liao2018sharp}, the approximation on the non-uniform mesh of the Volterra operator given in Eq.~\eqref{eq1.1} is
\begin{equation}
\begin{split}
\int_{0}^{t} \xi_{1 - \beta}(t - s) \partial_{s} u(x,y,s) ds \mid_{t = t_k} & \approx
a_{0}^{(k,\beta)} u(x,y,t_k)
+ \sum\limits_{\ell = 1}^{k - 1} \left( a_{k - \ell}^{(k,\beta)} - a_{k - \ell - 1}^{(k,\beta)} \right) u(x,y,t_{\ell}) \\
&\quad - a_{k - 1}^{(k,\beta)} u(x,y,t_{0}) \\
& \triangleq \delta_t^\beta u(x,y,t_k),
\end{split}
\label{eq2.1}
\end{equation}
where
\begin{equation*}
a_{k - \ell}^{(k,\beta)} = \int_{t_{\ell - 1}}^{t_\ell} \frac{\xi_{1 - \beta}(t_k - s)}{\tau_\ell} ds
= \frac{\xi_{2 - \beta}(t_k - t_{\ell - 1}) - \xi_{2 - \beta}(t_k - t_{\ell})}{\tau_\ell}, \quad \ell = 1, 2, \cdots, k.
\end{equation*}

Let $u_{i j}^{k}$ be the approximation of $u(x_i,y_j, t_k)$ and $f_{i j}^{k} = f(x_i, y_j, t_k)$.
Then, the time-stepping scheme of Eq.~\eqref{eq1.1} is
\begin{equation}
\delta_t^\beta u_{ij}^{k}
= \kappa \delta_{xy}^{2} u_{ij}^{k} + f_{ij}^{k} \quad\mathrm{for}\quad 
1 \leq i \leq N_x - 1,~1 \leq j \leq N_y - 1,~1 \leq k \leq M,
\label{eq2.2}
\end{equation}
where 
\begin{equation*}
\delta_{xy}^{2} u_{ij}^{k} = \frac{u_{i - 1,j}^{k} - 2 u_{ij}^{k} + u_{i + 1,j}^{k}}{h_x^2}
+ \frac{u_{i,j - 1}^{k} - 2 u_{ij}^{k} + u_{i,j + 1}^{k}}{h_y^2}.
\end{equation*}
Let 
\begin{equation*}
\begin{split}
& B_x = 1/h_x^2 \;\mathrm{tridiag} (1, -2, 1) \in \mathbb{R}^{(N_x - 1) \times (N_x - 1)}, \quad
B_y = 1/h_y^2 \;\mathrm{tridiag} (1, -2, 1) \in \mathbb{R}^{(N_y - 1) \times (N_y - 1)}, \\
& \bm{u}^k = \left[ u_{11}^{k}, \cdots, u_{N_x - 1, 1}^{k}, 
u_{12}^{k}, \cdots, u_{N_x - 1, 2}^{k}, \cdots, 
u_{1,N_y - 1}^{k}, \cdots, u_{N_x - 1,N_y - 1}^{k} \right]^{T}
\end{split}
\end{equation*}
and
\begin{equation*}
\bm{f}^k = \left[ f_{11}^{k}, \cdots, f_{N_x - 1, 1}^{k}, 
f_{12}^{k}, \cdots, f_{N_x - 1, 2}^{k}, \cdots, 
f_{1,N_y - 1}^{k}, \cdots, f_{N_x - 1,N_y - 1}^{k} \right]^{T}.
\end{equation*}
Then, the matrix form of Eq.~\eqref{eq2.2} is given by
\begin{equation}
\delta_t^\beta \bm{u}^k = B \bm{u}^k + \bm{f}^k,
\label{eq2.3}
\end{equation}
where $B = \kappa \left( I_y \otimes B_x + B_y \otimes I_x \right)$ (``$\otimes$" denotes the Kronecker product).
Here $I_x$ and $I_y$ are two identity matrices with sizes $N_x - 1$ and $N_y - 1$, respectively.
%$1/h_x^2 \mathrm{tridiag} (1, -2, 1)
Furthermore, the stability and convergence of the time-stepping scheme (\ref{eq2.2}) can be proved simply based on the work \cite{stynes2017error,liao2018sharp}.

\subsection{The all-at-once system}
\label{sec2.2}

Before deriving our all-at-once system,
several auxiliary symbols are introduced:
$I_t$ and $I_s$ represent identity matrices with sizes $M$ and $(N_x - 1)(N_y - 1)$, respectively.
Denote
\begin{equation*}
\bm{u} = \left[ \left( \bm{u}^1 \right)^T,  \left( \bm{u}^2 \right)^T,\cdots, \left( \bm{u}^M \right)^T \right]^T \quad\mathrm{and}\quad
\bm{f} = \left[ \left( \bm{f}^1 \right)^T,  \left( \bm{f}^2 \right)^T,\cdots, \left( \bm{f}^M \right)^T \right]^T.
\end{equation*}

With the help of Eq. (\ref{eq2.3}), the all-at-once system can be written as:
\begin{equation}{}
\mathcal{M} \bm{u} = \bm{\eta} + \bm{f},
\label{eq2.4}
\end{equation}
where $\mathcal{M} = A \otimes I_s - I_t \otimes B$ with
\begin{equation*}
A =
\begin{bmatrix}
a_0^{(1,\beta)} & 0 & 0 & \cdots & 0 & 0 \\
a_1^{(2,\beta)} - a_0^{(2,\beta)} & a_0^{(2,\beta)} & 0 & \cdots & \cdots & 0 \\
\vdots & a_1^{(3,\beta)} - a_0^{(3,\beta)} & a_0^{(3,\beta)} & \ddots & \ddots & \vdots \\
\vdots & \ddots & \ddots & \ddots & \ddots & 0 \\
a_{M - 2}^{(M - 1,\beta)} - a_{M - 3}^{(M - 1,\beta)} & \ddots & \ddots & \ddots & a_0^{(M - 1,\beta)} & 0 \\
a_{M - 1}^{(M,\beta)} - a_{M - 2}^{(M,\beta)} & a_{M - 2}^{(M,\beta)} - a_{M - 3}^{(M,\beta)} &
\cdots & \cdots &  a_1^{(M,\beta)} - a_0^{(M,\beta)} & a_0^{(M,\beta)}
\end{bmatrix}
\end{equation*}
and $\bm{\eta} = \left[ a_0^{(1,\beta)} \left( \bm{u}^0 \right)^T, a_1^{(2,\beta)} \left( \bm{u}^0 \right)^T, \cdots,
a_{M - 1}^{(M,\beta)} \left( \bm{u}^0 \right)^T\right]^T$.

If Gaussian elimination is applied
in the block forward substitution (BFS) method \cite{ke2015fast,huang2017fast} to solve Eq.~\eqref{eq2.4},
the matrix $\mathcal{M}$ must be stored.
Thus, the computational complexity and storage requirement of such a method are
$\mathcal{O}(M N_x^3 N_y^3 + M^2 N_x N_y)$ and $\mathcal{O}(M N_x^2 N_y^2)$, respectively.
%$\mathcal{O}(M N^2 + M^2 N)$ and $\mathcal{O}(M N^3 + M^2 N^2)$.
In order to reduce the computational cost,
we prefer to use Krylov subspace methods 
such as the biconjugate gradient stabilized (BiCGSTAB) method \cite{van1992bi}.
In the next section, an efficient algorithm is designed for fast solving Eq.~\eqref{eq2.4}.

\section{Two preconditioners and their spectral analysis}
\label{sec3}

In this section, two preconditioners are designed for solving Eq.~\eqref{eq2.4}.
Several properties of these preconditioners are also investigated.
It is easy to find that one part of $A$ has Toeplitz structure due to the uniform mesh is used in $[T_0,T]$.
Thus, the matrix $A$ can be rewritten as the following $2 \times 2$ block matrix:
\begin{align*}
A =
\begin{bmatrix}
A_{11} & \bm{0} \\
A_{21} & A_{22}
\end{bmatrix},
\end{align*}
where $\bm{0}$ is a zero matrix with suitable size,
\begin{align*}
A_{11} =
\begin{bmatrix}
a_0^{(1,\beta)} & 0 & 0 & \cdots & 0 & 0 \\
a_1^{(2,\beta)} - a_0^{(2,\beta)} & a_0^{(2,\beta)} & 0 & \cdots & \cdots & 0 \\
\vdots & a_1^{(3,\beta)} - a_0^{(3,\beta)} & a_0^{(3,\beta)} & \ddots & \ddots & \vdots \\
\vdots & \ddots & \ddots & \ddots & \ddots & 0 \\
a_{M_0 - 2}^{(M_0 - 1,\beta)} - a_{M_0 - 3}^{(M_0 - 1,\beta)} & \ddots & \ddots & \ddots & a_0^{(M_0 - 1,\beta)} & 0 \\
a_{M_0 - 1}^{(M_0,\beta)} - a_{M_0 - 2}^{(M_0,\beta)} & a_{M_0 - 2}^{(M_0,\beta)} - a_{M_0 - 3}^{(M_0,\beta)} &
\cdots & \cdots &  a_1^{(M_0,\beta)} - a_0^{(M_0,\beta)} & a_0^{(M_0,\beta)}
\end{bmatrix},
\end{align*}
\begin{equation*}
A_{21} =
\begin{bmatrix}
a_{M_0}^{(M_0 + 1,\beta)} - a_{M_0 - 1}^{(M_0 + 1,\beta)} & a_{M_0 - 1}^{(M_0 + 1,\beta)} - a_{M_0 - 2}^{(M_0 + 1,\beta)} & \cdots
& a_{1}^{(M_0 + 1,\beta)} - b_{0}^{(\beta)}\\
a_{M_0 + 1}^{(M_0 + 2,\beta)} - a_{M_0}^{(M_0 + 2,\beta)} & a_{M_0}^{(M_0 + 2,\beta)} - a_{M_0 - 1}^{(M_0 + 2,\beta)} & \cdots
& a_{2}^{(M_0 + 2,\beta)} - b_{1}^{(\beta)} \\
\vdots & \vdots &  & \vdots \\
a_{M - 2}^{(M - 1,\beta)} - a_{M - 3}^{(M - 1,\beta)} & a_{M - 3}^{(M - 1,\beta)} - a_{M - 4}^{(M - 1,\beta)} & \cdots
&a_{M - M_0 - 1}^{(M - 1,\beta)} - b_{M - M_0 - 2}^{(\beta)} \\
a_{M - 1}^{(M,\beta)} - a_{M - 2}^{(M,\beta)} & a_{M - 2}^{(M,\beta)} - a_{M - 3}^{(M,\beta)} & \cdots
&a_{M - M_0}^{(M,\beta)} - b_{M - M_0 - 1}^{(\beta)}
\end{bmatrix}
\end{equation*}
and
\begin{equation*}
A_{22} =
\begin{bmatrix}
\omega_{0}^{(\beta)} & 0 & 0 & \cdots & 0 & 0 \\
\omega_{1}^{(\beta)} & \omega_{0}^{(\beta)} &  0 & \cdots & \cdots & 0 \\
\vdots & \omega_{1}^{(\beta)} & \omega_{0}^{(\beta)} & \ddots & \ddots & \vdots \\
\vdots & \ddots & \ddots & \ddots & \ddots & 0 \\
\omega_{M - M_0 - 2}^{(\beta)} & \ddots & \ddots & \ddots & \omega_0^{(\beta)} & 0 \\
\omega_{M - M_0 - 1}^{(\beta)} & \omega_{M - M_0 - 2}^{(\beta)} & \cdots & \cdots & \omega_1^{(\beta)} & \omega_0^{(\beta)}
\end{bmatrix}
\end{equation*}
with
\begin{equation}
\omega_{k}^{(\beta)} =
\begin{cases}
b_0^{(\beta)}, & k = 0, \\
b_{k}^{(\beta)} - b_{k - 1}^{(\beta)}, & k = 1,2,\cdots, M - M_0 - 1.
\end{cases}
\label{eq3.1}
\end{equation}
Here, $b_{\ell}^{(\beta)} = \frac{\tilde{\tau}^{-\beta}}{\Gamma(2 - \beta)} \left[ (\ell + 1)^{1 - \beta} - \ell^{1 - \beta} \right]$ for $\ell \geq 0$.
Let $I_{t1}$ and $I_{t2}$ be two identity matrices with sizes $M_0$ and $M - M_0$, respectively.
Then, the solution of Eq.~\eqref{eq2.4} can be obtained by solving the following equivalent two subproblems:
\begin{subequations}
\begin{align}
\mathcal{M}_{11} \tilde{\bm{u}}_1 &= \bm{\eta}_1 + \tilde{\bm{f}}_1, \label{eq3.2a} \\
\mathcal{M}_{22} \tilde{\bm{u}}_2 &= \bm{\eta}_2 + \tilde{\bm{f}}_2 - \mathcal{M}_{21} \tilde{\bm{u}}_1,
\label{eq3.2b}
\end{align}
\label{eq3.2}
\end{subequations}
where
%$\bm{U}^1 = \left[ \left( \bm{u}^1 \right)^T,  \left( \bm{u}^2 \right)^T,\cdots, \left( \bm{u}^{M_0} \right)^T \right]^T$,
%$\bm{U}^2 = \left[ \left( \bm{u}^{M_0 + 1} \right)^T,  \left( \bm{u}^{M_0 + 2} \right)^T,\cdots, \left( \bm{u}^{M} \right)^T \right]^T$,
\begin{align*}
& \mathcal{M}_{11} = A_{11} \otimes I_s - I_{t1} \otimes B, \quad \mathcal{M}_{21} = A_{21} \otimes I_s,
\quad \mathcal{M}_{22} = A_{22} \otimes I_s - I_{t2} \otimes B, \\
& \tilde{\bm{u}}_1 = \left[ \left( \bm{u}^1 \right)^T,  \left( \bm{u}^2 \right)^T,\cdots, \left( \bm{u}^{M_0} \right)^T \right]^T, \quad
\tilde{\bm{u}}_2 = \left[ \left( \bm{u}^{M_0 + 1} \right)^T,  \left( \bm{u}^{M_0 + 2} \right)^T,\cdots, \left( \bm{u}^{M} \right)^T \right]^T, \\
& \bm{\eta}_1 = \left[ a_0^{(1,\beta)} \left( \bm{u}^0 \right)^T, a_1^{(2,\beta)} \left( \bm{u}^0 \right)^T, \cdots,
a_{M_0 - 1}^{(M_0,\beta)} \left( \bm{u}^0 \right)^T \right]^T, \\
& \bm{\eta}_2 = \left[ a_{M_0}^{(M_0 + 1,\beta)} \left( \bm{u}^0 \right)^T,
a_{M_0 + 1}^{(M_0 + 2,\beta)} \left( \bm{u}^0 \right)^T, \cdots,
a_{M - 1}^{(M,\beta)} \left( \bm{u}^0 \right)^T \right]^T, \\
& \tilde{\bm{f}}_1 = \left[ \left( \bm{f}^1 \right)^T,\left( \bm{f}^2 \right)^T,\cdots,\left( \bm{f}^{M_0} \right)^T \right]^T, \quad
\tilde{\bm{f}}_2 = \left[ \left( \bm{f}^{M_0 + 1} \right)^T,\left( \bm{f}^{M_0 + 2} \right)^T,\cdots,\left( \bm{f}^{M} \right)^T \right]^T.
\end{align*}

These two subproblems are only coupled by $\tilde{\bm{u}}_1$.
Thus, different methods can be used to solve them, such as
the DC method \cite{ke2015fast,huang2017fast} and the AI method \cite{lu2015fast,lu2018approximate}.
For example, Krylov subspace methods \cite{saad2003iterative} are used to solve Eq.~\eqref{eq3.2a},
while the DC method \cite{ke2015fast,huang2017fast} is employed to solve Eq.~\eqref{eq3.2b}.
In this work, we choose the preconditioned BiCGSTAB (PBiCGSTAB) method \cite{van1992bi} to solve both of them.
It is worth mentioning that the general minimum residual (GMRES) method \cite{saad2003iterative} 
	is not used, since it requires large amounts of storage due to the orthogonalization process.

\subsection{A block lower tridiagonal preconditioner for Eq.~(3.2a)}
\label{sec3.1}

\begin{figure}[t]
	\centering
	\subfigure[$(r, M_0) = (2,7)$]
	{\includegraphics[width=2.6in,height=2.5in]{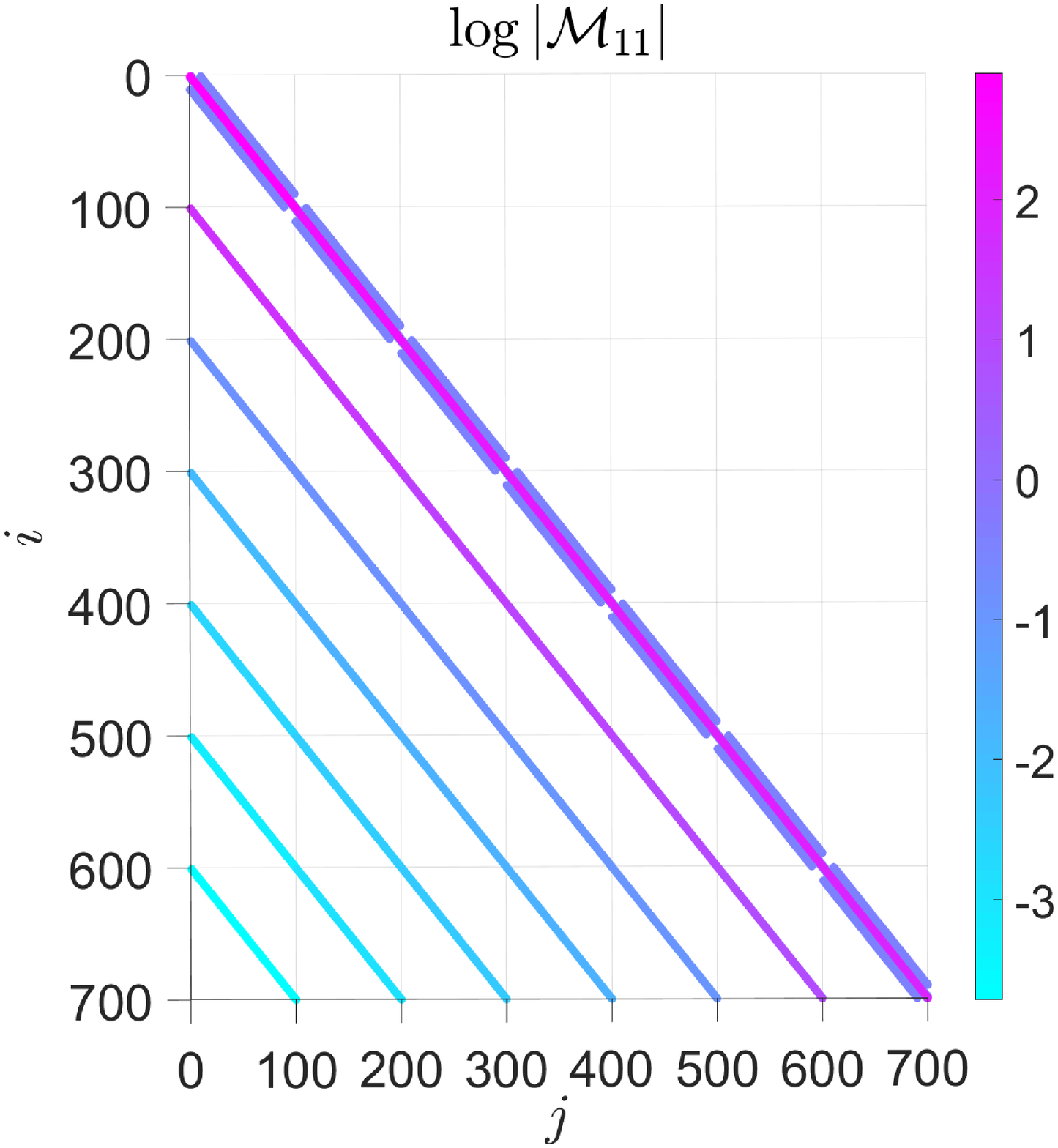}}\hspace{6mm}
	%\caption{Sparsity pattern of $W \in \mathbb{R}^{100 \times 100}$, when $M = N = 11$.}
	\subfigure[$(r, M_0) = (3,5)$]
	{\includegraphics[width=2.5in,height=2.5in]{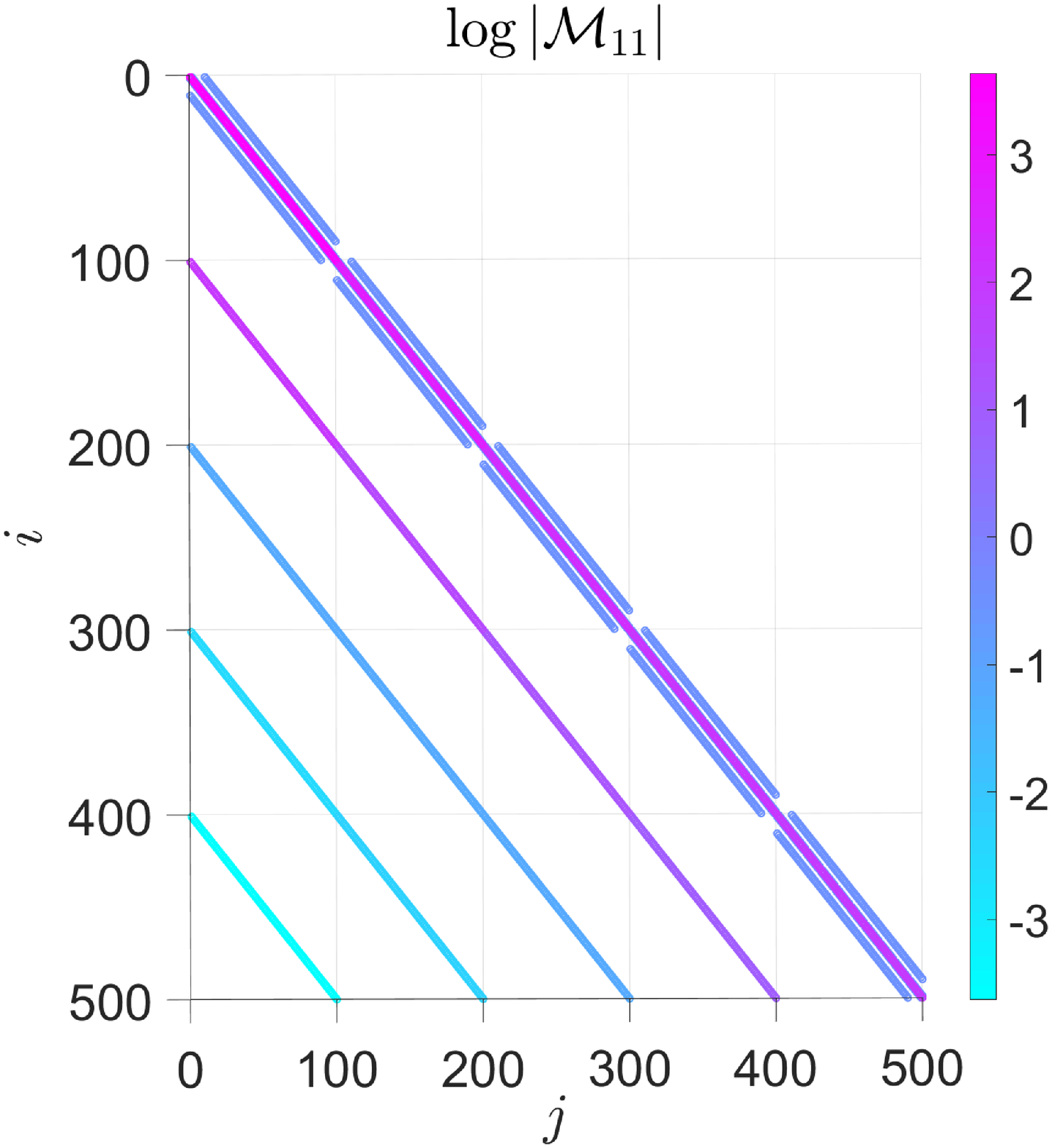}}
	\caption{The decay of the elements of matrix $\mathcal{M}_{11}$, where $\beta = 0.5$ and $N_x = N_y = 11$.}
	\label{fig1}
\end{figure}
It can be seen from Fig.~\ref{fig1} that the diagonal entries of $\mathcal{M}_{11}$ decay quickly. 
Inspired by this observation, a block lower tridiagonal preconditioner is designed for solving the subproblem \eqref{eq3.2a}:
\begin{equation*}
P_1 = \mathrm{tri}(A_{11}) \otimes I_s - I_{t1} \otimes B,
\end{equation*}
where $\mathrm{tri}(A_{11})$ is a matrix which only preserves the first three diagonals of $A_{11}$ 
and the others are set $0$, i.e.,
\begin{equation*}
\text{tri}(A_{11}) =
\begin{bmatrix}
a_0^{(1,\beta)} & 0 & 0 & \cdots & 0 & 0 \\
a_1^{(2,\beta)} - a_0^{(2,\beta)} & a_0^{(2,\beta)} & 0 & \cdots & \cdots & 0 \\
a_2^{(3,\beta)} - a_1^{(3,\beta)} & a_1^{(3,\beta)} - a_0^{(3,\beta)} & a_0^{(3,\beta)} & \ddots & \ddots & \vdots \\
0 & \ddots & \ddots & \ddots & \ddots & 0 \\
\vdots & \ddots & \ddots & \ddots & a_0^{(M_0 - 1,\beta)} & 0 \\
0 & 0 & \cdots & \cdots &  a_1^{(M_0,\beta)} - a_0^{(M_0,\beta)} & a_0^{(M_0,\beta)}
\end{bmatrix}.
\end{equation*}

The nonsingularity of $P_1$ is easy to check since all main diagonal blocks of it are nonsingular.

\begin{theorem}
The degree of the minimal polynomial $p$ \cite{murphy2000note} of $P_1^{-1} \mathcal{M}_{11}$ satisfies that
\begin{equation*}
\mathrm{deg} \; p(P_1^{-1} \mathcal{M}_{11}) \leq \lceil M_0/3 \rceil.
%= \left( P_1^{-1} \mathcal{M}_{11} - I_{ts1} \right)^{\lceil M_0/2\rceil},
\end{equation*}
%where $I_{ts1} = I_{t1} \otimes I_s$.
Thus, the dimension of the Krylov subspace 
	$\mathcal{K} \left( P_1^{-1} \mathcal{M}_{11}; \bm{b} \right)$ is at most $\lceil M_0/3 \rceil$.
\label{th3.1}
\end{theorem}
% Proof
\textbf{Proof.}
After some simple calculations, we have
\begin{equation*}
P_1^{-1} \mathcal{M}_{11} =
\begin{bmatrix}
I_{s} & \bm{0} & \bm{0} & \cdots & \cdots & \cdots & \cdots & \bm{0} \\
\bm{0} & I_{s} & \bm{0} & \cdots & \cdots & \cdots & \cdots & \vdots\\
\bm{0} & \bm{0} & I_{s} & \ddots &  &  &  & \vdots \\
J_{4}^{4} & \bm{0} & \ddots & \ddots & \ddots &  &  & \vdots\\
J_{5}^{5} & J_{5}^{4} & \ddots & \ddots & \ddots & \ddots &  & \vdots\\
\vdots & \ddots & \ddots & \ddots & \ddots & \ddots & \ddots & \vdots \\
J_{M_0 - 1}^{M_0 - 1} & \ddots & \ddots & \ddots & \ddots & \ddots & \ddots & \bm{0} \\
J_{M_0}^{M_0} & J_{M_0}^{M_0 - 1} & \cdots & \cdots & J_{M_0}^{4} & \bm{0} & \bm{0} & I_{s}
\end{bmatrix},
\end{equation*}
where
%for $m \leq k \leq M_0$,
%\begin{equation*}
%J_{k}^{s} = 	
%\begin{cases}
%\left( a_{3}^{(k,\beta)} - a_{2}^{(k,\beta)} \right) \left(a_{0}^{(k,\beta)} I_s - B \right)^{-1}, & m = 4, \\
%\left(a_{0}^{(k,\beta)} I_s - B \right)^{-1} \left[ \left( a_{4}^{(k,\beta)} - a_{3}^{(k,\beta)} \right) I_s
%- \left( a_{1}^{(k,\beta)} - a_{0}^{(k,\beta)} \right) J_{k - 1}^{4} \right], & m = 5, \\
%\left(a_{0}^{(k,\beta)} I_s - B \right)^{-1} \left[ \left( a_{s - 1}^{(k,\beta)} - a_{s - 2}^{(k,\beta)} \right) I_s
%- \left( a_{2}^{(k,\beta)} - a_{1}^{(k,\beta)} \right) J_{k - 2}^{s - 2}
%- \left( a_{1}^{(k,\beta)} - a_{0}^{(k,\beta)} \right) J_{k - 1}^{s - 1} \right], & 6 \leq m \leq M_0.	
%\end{cases}
%\end{equation*}
%
\begin{align*}
& J_{k}^{4} = \left( a_{3}^{(k,\beta)} - a_{2}^{(k,\beta)} \right) \left(a_{0}^{(k,\beta)} I_s - B \right)^{-1}, \quad k \geq 4, \\
& J_{k}^{5} = \left(a_{0}^{(k,\beta)} I_s - B \right)^{-1} \left[ \left( a_{4}^{(k,\beta)} - a_{3}^{(k,\beta)} \right) I_s
- \left( a_{1}^{(k,\beta)} - a_{0}^{(k,\beta)} \right) J_{k - 1}^{4} \right], \quad k \geq 5,
\end{align*}
and, for $6 \leq m \leq M_0$,
\begin{align*}
J_{k}^{m} = \left(a_{0}^{(k,\beta)} I_s - B \right)^{-1} \left[ \left( a_{k - 1}^{(k,\beta)} - a_{k - 2}^{(k,\beta)} \right) I_s
- \left( a_{2}^{(k,\beta)} - a_{1}^{(k,\beta)} \right) J_{k - 2}^{m - 2}
- \left( a_{1}^{(k,\beta)} - a_{0}^{(k,\beta)} \right) J_{k - 1}^{m - 1} \right], ~ m \leq k \leq M_0.
\end{align*}

This implies that all eigenvalues of $P_1^{-1} \mathcal{M}_{11}$ are equal to $1$.
On the other hand, $P_1^{-1} \mathcal{M}_{11} - I_{ts1}$ is a strictly block lower triangular matrix, where $I_{ts1} = I_{t1} \otimes I_s$.
A simple computation shows that $\left( P_1^{-1} \mathcal{M}_{11} - I_{ts1} \right)^{\lceil M_0/3 \rceil} = \bm{0}$.
However, it is even possible that $\left( P_1^{-1} \mathcal{M}_{11} - I_{ts1} \right)^{m} = \bm{0}$ 
	for some $m \leq \lceil M_0/3 \rceil$.
This means that the minimal polynomial $p(P_1^{-1} \mathcal{M}_{11})$ has a maximum degree of $\lceil M_0/3 \rceil$.
From Saad's work \cite{saad2003iterative}, we know that the dimension of the Krylov
subspace $\mathcal{K} \left( P_1^{-1} \mathcal{M}_{11}; \bm{b} \right)$ is then also at most $\lceil M_0/3 \rceil$.
Thus, the proof is completed.
\hfill $\Box$

According to Theorem \ref{th3.1}, it can be seen that if a Krylov subspace method 
	with an optimal or Galerkin property is employed to solve the preconditioned form of \eqref{eq3.2a}
with the coefficient matrix $P_1^{-1} \mathcal{M}_{11}$ in exact arithmetic, 
it will converge to the exact solution of \eqref{eq3.2a}
in at most $\lceil M_0/3 \rceil$ iterations. 

In practice, preconditioned Krylov subspace methods need to calculate $P_1^{-1} \bm{v}$, 
where $\bm{v}$ is a vector.
Let 
\begin{equation*}
\begin{split}
& Q_x = \left( \sqrt{2/N_x} \sin \left( \frac{i j \pi}{N_x} \right) \right)_{1 \leq i,j \leq N_x - 1}, \qquad
Q_y = \left( \sqrt{2/N_y} \sin \left( \frac{i j \pi}{N_y} \right) \right)_{1 \leq i,j \leq N_y - 1}, \\
& D_{Bx} = \mathrm{diag} \left( \lambda_{1}^{Bx},\lambda_{2}^{Bx}, \cdots,\lambda_{N_x - 1}^{Bx} \right), \qquad
D_{By} = \mathrm{diag} \left( \lambda_{1}^{By},\lambda_{2}^{By}, \cdots,\lambda_{N_y - 1}^{By} \right)
\end{split}
\end{equation*}
with
\begin{equation*}
\lambda_{i}^{Bx} = -\frac{4 \kappa}{h_x^2} \sin^2 \left( \frac{i \pi}{2 N_x} \right) < 0 \quad \mathrm{for} \quad
1 \leq i \leq N_x - 1
\end{equation*}
and 
\begin{equation*}
\lambda_{j}^{By} = -\frac{4 \kappa}{h_y^2} \sin^2 \left( \frac{j \pi}{2 N_y} \right) < 0 \quad \mathrm{for} \quad
1 \leq j \leq N_y - 1.
\end{equation*}
According to \cite[Sec.~4.3]{moroney2013efficient}, we know that
$B = Q D_{B} Q^T$, where
\begin{equation*}
Q = Q_y \otimes Q_x \quad \mathrm{and} \quad
D_{B} = I_y \otimes D_{Bx} + D_{By} \otimes I_x
= \mathrm{diag} \left( \lambda_{1}^{B},\lambda_{2}^{B}, \cdots,\lambda_{(N_x - 1) (N_y - 1)}^{B} \right).
\end{equation*}
Hence, $\bm{z} = P_{1}^{-1} \bm{v}$ can be computed in the following way:
\begin{equation}
\begin{cases}
\tilde{\bm{z}}_1 = \left( I_{t1} \otimes Q \right) \bm{v}, & \textrm{Step-(a1)}, \\
\left[ \mathrm{tri}(A_{11}) \otimes I_s - I_{t1} \otimes D_{B}  \right] \tilde{\bm{z}}_{2} = \tilde{\bm{z}}_{1}, & \textrm{Step-(b1)}, \\
\bm{z} = \left(  I_{t1} \otimes Q^T \right)  \tilde{\bm{z}}_2, & \textrm{Step-(c1)}.
\end{cases}
\label{eq3.51}
\end{equation}

In Eq.~\eqref{eq3.51}, the first and the third step can be done by discrete sine transform.
Each step itself can be carried out in parallel on $M_0$ processors.
We notice that the matrix in $\textrm{Step-(b1)}$
is a block tridiagonal matrix with diagonal blocks. 
Thus, $\tilde{\bm{z}}_{2}$ can be obtained by the BFS method 
with $\mathcal{O}(M_0 N_x N_y)$ operations and $\mathcal{O}(M_0 N_x N_y)$ memory.
%The banded $LU$ factorization \cite{demmel1997applied}: $P_1 = L_1 U_1$ is used in our paper.
%Then, $P_1^{-1} \bm{v} = U_1^{-1} \left( L_1^{-1} \bm{v} \right)$.
%The full cost of computing $P_1^{-1} \bm{v}$ is much less than $\mathcal{O}(M_0 N^2 + N)$ since $P_1$ is very sparse.
%Moreover, we choose the MATLAB software package CSparse~\footnote{http://faculty.cse.tamu.edu/davis/SuiteSparse/} for fast computing $P_1^{-1} \bm{v}$.

\subsection{A PinT preconditioner for Eq.~(3.2b)}
\label{sec3.2}

In this subsection, we concentrate on solving Eq.~\eqref{eq3.2b} efficiently.
Note that this subproblem is a block triangular Toeplitz matrix with a tridiagonal block system.
The DC method \cite{ke2015fast,huang2017fast} or the AI method \cite{lu2015fast,lu2018approximate} can be used to solve it.
In this paper, the system is solved by a Krylov subspace method with a preconditioner.
This preconditioner can be efficiently implemented in the PinT framework.
Our PinT preconditioner (or structuring-circulant preconditioner, or semi-circulant preconditioner) is
\begin{equation}
P_{\alpha} = A_{22}^{\alpha} \otimes I_s - I_{t2} \otimes B,
\label{eq3.3}
\end{equation}
where $A_{22}^{\alpha} = A_{22} + \alpha \tilde{A}$ is a $\alpha$-circulant matrix with the parameter $\alpha \in (0,1]$ and
\begin{equation*}
\tilde{A} =
\begin{bmatrix}
0 & \omega_{M - M_0 - 1}^{(\beta)} & \cdots & \omega_2^{(\beta)} & \omega_1^{(\beta)} \\
0 & 0 & \omega_{M - M_0 - 1}^{(\beta)} & \cdots & \omega_2^{(\beta)} \\
0 & 0 & 0 & \ddots & \vdots \\
\vdots & \ddots & \ddots & \ddots & \omega_{M - M_0 - 1}^{(\beta)} \\
0 & \cdots & 0 & 0 & 0
\end{bmatrix}.
\end{equation*}

Denote $\Theta_{\alpha} = \mathrm{diag} \left( 1, \alpha^{-\frac{1}{M - M_0}}, \cdots, \alpha^{-\frac{M - M_0 - 1}{M - M_0 }} \right)$.
Let $\mathbb{F}$ represent the discrete Fourier matrix, let ``$^*$" denote the conjugate transpose of a matrix.
We know that the $\alpha$-circulant matrix $A_{22}^{\alpha}$ has the following diagonalization:
\begin{equation*}
A_{22}^{\alpha} = V_{\alpha} \Lambda_{\alpha} V_{\alpha}^{-1}
\end{equation*}
with $V_{\alpha} = \Theta_{\alpha} \mathbb{F}^{*}$ and
\begin{equation}
\Lambda_{\alpha} = \mathrm{diag} \left( \mathbb{F} \Theta_{\alpha}^{-1} A_{22}^{\alpha}(:,1) \right)
= \mathrm{diag} \left( \lambda_{1}^{(\alpha)},\lambda_{2}^{(\alpha)}, \cdots,\lambda_{M - M_0}^{(\alpha)} \right)
\label{eq3.4}
\end{equation}
contains all eigenvalues of $A_{22}^{\alpha}$, where $A_{22}^{\alpha}(:,1)$ is the first column of $A_{22}^{\alpha}$.

The diagonalization of $A_{22}^{\alpha}$ immediately implies that
$\bm{z} = P_{\alpha}^{-1} \bm{v}$ can be computed via the following three steps:
\begin{equation}
\begin{cases}
\bm{z}_1 = \left( V_\alpha^{-1} \otimes I_s \right) \bm{v}
= \left( \mathbb{F} \otimes I_s \right) \left[ \left( \Theta_{\alpha}^{-1} \otimes I_s \right) \bm{v} \right], & \textrm{Step-(a)}, \\
\left( \lambda_n^{(\alpha)} I_s - B  \right) \bm{z}_{2,n} = \bm{z}_{1,n}, ~1 \leq n \leq M - M_0, & \textrm{Step-(b)}, \\
\bm{z} = \left( V_\alpha \otimes I_s \right) \bm{z}_2
= \left(  \Theta_{\alpha} \otimes I_s \right) \left[ \left( \mathbb{F}^{*} \otimes I_s \right) \bm{z}_2 \right], & \textrm{Step-(c)},
\end{cases}
\label{eq3.5}
\end{equation}
where $\bm{z}_j = \left[ \bm{z}_{j,1}^T, \bm{z}_{j,2}^T, \cdots, \bm{z}_{j,M - M_0}^T \right]^T$ with $j = 1,2$.
In the first and third step, the matrix-vector multiplications can be done
by fast Fourier transforms (FFTs)~\footnote{A parallel version of FFT is available at \url{http://www.fftw.org/parallel/parallel-fftw.html}.}
in parallel on $M - M_0$ processors \cite{gu2020parallel}.
Thus, the major computation cost of $P_{\alpha}^{-1} \bm{v}$ comes from the second step.
In this step, $M - M_0$ complex linear equations need to be solved, but they can be solved simultaneously (i.e., parallel computing).
\begin{remark}
According to the work \cite{gu2021pint}, we only solve the first $\left \lceil \frac{M - M_0 + 1}{2} \right \rceil$
shifted complex-valued linear systems in $\textrm{Step-(b)}$.
On the other hand, the matrix $B$ can be diagonalized by a discrete sine transform.
Thus, the storage requirement and the computational cost in $\textrm{Step-(b)}$ are
$\mathcal{O} \left(\left \lceil \frac{M - M_0 + 1}{2} \right \rceil N_x N_y \right)$
and $\mathcal{O} \left( \left \lceil \frac{M - M_0 + 1}{2} \right \rceil N_x N_y \log N_x N_y \right)$, respectively.
\end{remark}

Firstly, we investigate the nonsingularity of $P_{\alpha}$. For this, the following result is needed.
\begin{lemma}
For any $\beta \in (0,1)$ and $\omega_{k}^{(\beta)}$ defined in \eqref{eq3.1}, it holds
\begin{equation}
\omega_{0}^{(\beta)} > \sum_{\ell = 1}^{M - M_0 - 1} \left| \omega_{\ell}^{(\beta)} \right| \quad \mathrm{and} \quad
\omega_{\ell}^{(\beta)} < 0 \quad \mathrm{for} \quad \ell \geq 1.
\label{eq3.6}
\end{equation}
\label{lemma3.1}
\end{lemma}
% Proof
\textbf{Proof.}
It is direct to check that
\begin{equation*}
\frac{\tilde{\tau}^{-\beta}}{\Gamma(2 - \beta)} = \omega_{0}^{(\beta)} = b_{0}^{(\beta)} > b_{1}^{(\beta)} > \cdots > b_{\ell}^{(\beta)} > 0, \quad
b_{\ell}^{(\beta)} \rightarrow 0 \quad \mathrm{as} \quad \ell \rightarrow +\infty.
\end{equation*}
Immediately, we have $\omega_{\ell}^{(\beta)} < 0$ for $\ell \geq 1$.
Then,
\begin{equation*}
\sum_{\ell = 1}^{M - M_0 - 1} \left| \omega_{\ell}^{(\beta)} \right|
= \sum_{\ell = 1}^{M - M_0 - 1} \left( b_{\ell - 1}^{(\beta)} - b_{\ell}^{(\beta)} \right)
= \omega_{0}^{(\beta)} - b_{M - M_0 - 1}^{(\beta)} < \omega_{0}^{(\beta)},
\end{equation*}
which completes the proof.
$\hfill\Box$

Based on Eq.~\eqref{eq3.5} and the diagonalization of $B$, 
in order to prove the nonsingularity of $P_{\alpha}$, we only need to show that
the real part of $\lambda_{j}^{(\alpha)}$ is positive, i.e., $\text{Re\,}(\lambda_{j}^{(\alpha)}) > 0$ for $1 \leq j \leq M - M_0$.
\begin{theorem}
For any $\beta \in (0,1)$ and $\alpha \in (0,1]$, it holds that
\begin{equation*}
\text{Re\,}(\lambda_{j}^{(\alpha)}) > 0 \quad \mathrm{for} \quad 1 \leq j \leq M - M_0.
\end{equation*}
\label{th3.2}
\end{theorem}
% Proof
\textbf{Proof.}
By Eq.~\eqref{eq3.4},
\begin{equation*}
\lambda_{j}^{(\alpha)} = \sum_{n = 0}^{M - M_0 - 1} \theta^{(j - 1) n} \alpha^{n/(M - M_0)} \omega_{n}^{(\beta)}, \quad 1 \leq j \leq M - M_0
\end{equation*}
with $\theta = \exp(\frac{2 \pi \iota}{M - M_0})$ and $\iota = \sqrt{-1}$.
For some given suitable angle $\phi_{j,n}$, we rewrite $\theta^{(j - 1) n} = \cos(\phi_{j,n}) + \iota \sin(\phi_{j,n})$.
Then, using Lemma \ref{lemma3.1}, we get
\begin{equation*}
\begin{split}
\text{Re\,}(\lambda_{j}^{(\alpha)}) & = \omega_{0}^{(\beta)} + \sum_{n = 1}^{M - M_0 - 1} \cos(\phi_{j,n}) \alpha^{n/(M - M_0)} \omega_{n}^{(\beta)} \\
& \geq \omega_{0}^{(\beta)} - \sum_{n = 1}^{M - M_0 - 1} \left| \cos(\phi_{j,n}) \alpha^{n/(M - M_0)} \omega_{n}^{(\beta)} \right| \\
& \geq \omega_{0}^{(\beta)} - \sum_{\ell = 1}^{M - M_0 - 1} \left| \omega_{\ell}^{(\beta)} \right| > 0,
\end{split}
\end{equation*}
and the proof is completed.
$\hfill\Box$

To estimate $I_{ts2} - P_{\alpha}^{-1} \mathcal{M}_{22}$ (where $I_{ts2} = I_{t2} \otimes I_s$),
the following lemma is useful.
\begin{lemma}(\cite{gu2020parallel})
For a strictly diagonally dominant (SDD) matrix $W \in \mathbb{C}^{M \times M}$, it holds
\begin{equation*}
\left\| W^{-1} \right\|_{\infty} \leq \frac{\max\limits_{1 \leq n \leq M} \tilde{z}_n(W)/\left| W(n,n) \right|}
{\min\limits_{1 \leq n \leq M} \left( 1 - \tilde{h}_n(W)/\left| W(n,n) \right| \right)},
\end{equation*}
where
$\left\{ \tilde{z}_n(W) \right\}_{n = 1}^{M}$ and $\left\{ \tilde{h}_n(W) \right\}_{n = 1}^{M}$ are given by the following two recursions
\begin{equation*}
\tilde{z}_1(W) = 1, \quad \tilde{z}_n(W) = 1 + \sum_{k = 1}^{n - 1} \frac{\left| W(n,k) \right|}{\left| W(k,k) \right|} \tilde{z}_k(W)
\end{equation*}
and
\begin{equation*}
\tilde{h}_1(W) = \sum_{k = 2}^{M} \left| W(1,k) \right|, \quad
\tilde{h}_n(W) = \sum_{k = 1}^{n - 1} \frac{\left| W(n,k) \right|}{\left| W(k,k) \right|} \tilde{h}_k(W) + \sum_{k = n + 1}^{M} \left| W(n,k) \right|,
\end{equation*}
respectively.
\label{lemma3.2}
\end{lemma}

For $\tilde{W} \in \mathbb{C}^{(M - M_0) (N_x - 1) (N_y - 1) \times (M - M_0) (N_x - 1) (N_y - 1)}$, we define the norm:
\begin{equation*}
\left\| \tilde{W} \right\|_{Q,\infty} = \left\| \left(I_{t2} \otimes Q^T \right) \tilde{W} \left(I_{t2} \otimes Q \right) \right\|_{\infty}.
\end{equation*}
In order to get a sharp estimate of $\left\| I_{ts2} - P_{\alpha}^{-1} \mathcal{M}_{22} \right\|_{Q,\infty}$,
we also need an auxiliary matrix defined as:
\begin{equation*}
L_\epsilon =
\begin{bmatrix}
1 &  &  &  & \\
\epsilon & 1 &  &  & \\
\epsilon^2 & \epsilon & 1 &  & \\
\vdots & \ddots & \ddots & \ddots & \\
\epsilon^{M - M_0} & \cdots & \epsilon^2 & \epsilon & 1
\end{bmatrix},
\end{equation*}
where $\epsilon$ is a free parameter.
%It is easy to know $\mu < 0$.
With this at hand, we can prove that $P_{\alpha}^{-1} \mathcal{M}_{22}$ is close to the identity.
\begin{theorem}
Let $\epsilon_{max} = - \omega_1^{(\beta)} / \omega_0^{(\beta)}$ and $R_{\epsilon_{max}} = L_{\epsilon_{max}} A_{22}^{\alpha}$.
Then, for any $\beta \in (0,1)$ and $\alpha \in (0,1]$, the following inequality holds
\begin{equation*}
\left\| I_{ts2} - P_{\alpha}^{-1} \mathcal{M}_{22} \right\|_{Q,\infty}
\leq C \alpha,
\end{equation*}
where
\begin{equation*}
C = \frac{\max\limits_{1 \leq n \leq M - M_0} \tilde{z}_n(R_{\epsilon_{max}})/ R_{\epsilon_{max}}(n,n) }
{\min\limits_{1 \leq n \leq M - M_0} \left( 1 - \tilde{h}_n(R_{\epsilon_{max}})/ R_{\epsilon_{max}}(n,n) \right)}
\sum_{n = 1}^{M - M_0 - 1} \epsilon_{max}^n \sum_{k = M - M_0 - n}^{M - M_0 - 1} \left| \omega_k^{(\beta)} \right|
\end{equation*}
is independent of the eigenvalues of $B$.
\label{th3.3}
\end{theorem}
% Proof
\textbf{Proof.}
From the definition of $P_\alpha$ \eqref{eq3.3}, we have
\begin{equation*}
\left\| I_{ts2} - P_{\alpha}^{-1} \mathcal{M}_{22} \right\|_{Q,\infty}
= \alpha \left\| \left(A_{22}^{\alpha} \otimes I_s  - I_{t2} \otimes D_{B} \right)^{-1} \left(\tilde{A} \otimes I_s \right) \right\|_{\infty}
= \alpha \max\limits_{\mu \in \sigma(B)} \left\| \left(A_{22}^{\alpha}  - \mu I_{t2} \right)^{-1} \tilde{A} \right\|_{\infty},
\end{equation*}
where $\sigma(B) = \left\{ \lambda_{k}^{B} \right\}_{k = 1}^{(N_x - 1)(N_y - 1)}$.
Thus, it turns to estimate
\begin{equation*}
\left\| \left(A_{22}^{\alpha}  - \mu I_{t2} \right)^{-1} \tilde{A} \right\|_{\infty}
= \left\| \left[ L_\epsilon \left(A_{22}^{\alpha}  - \mu I_{t2} \right) \right]^{-1} L_\epsilon \tilde{A} \right\|_{\infty}
\leq \left\| \left[ L_\epsilon \left(A_{22}^{\alpha}  - \mu I_{t2} \right) \right]^{-1} \right\|_{\infty} \left\| L_\epsilon \tilde{A} \right\|_{\infty}.
\end{equation*}
Denote $\tilde{R}_\epsilon = L_\epsilon \left(A_{22}^{\alpha}  - \mu I_{t2} \right)$,
$\tilde{A}_\epsilon = L_\epsilon \tilde{A}$
and let $\epsilon = - \frac{\omega_1^{(\beta)}}{\omega_0^{(\beta)} - \mu}$.
Then, according to Lemma 2.2 in \cite{gu2020parallel}, we know that
$\tilde{R}_\epsilon$ is a SDD matrix. Using Lemma \ref{lemma3.2}, we obtain
\begin{equation*}
\left\| \tilde{R}_\epsilon^{-1} \right\|_{\infty}
\leq \frac{\max\limits_{1 \leq n \leq M - M_0} \tilde{z}_n(\tilde{R}_\epsilon)/\left| \tilde{R}_\epsilon(n,n) \right|}
{\min\limits_{1 \leq n \leq M - M_0} \left( 1 - \tilde{h}_n(\tilde{R}_\epsilon)/\left| \tilde{R}_\epsilon(n,n) \right| \right)}
\leq \frac{\max\limits_{1 \leq n \leq M - M_0} \tilde{z}_n(R_{\epsilon_{max}})/ R_{\epsilon_{max}}(n,n) }
{\min\limits_{1 \leq n \leq M - M_0} \left( 1 - \tilde{h}_n(R_{\epsilon_{max}})/ R_{\epsilon_{max}}(n,n) \right)},
\end{equation*}
where the relation $\left| \tilde{R}_\epsilon(n,n) \right| \geq R_{\epsilon_{max}}(n,n)$ is used in the second inequality.

On the other hand, after a routine calculation, we get
\begin{equation*}
\left\| \tilde{A}_\epsilon \right\|_{\infty}
= \sum_{n = 1}^{M - M_0 - 1} \epsilon^n \sum_{k = M - M_0 - n}^{M - M_0 - 1} \left| \omega_k^{(\beta)} \right|
\leq \sum_{n = 1}^{M - M_0 - 1} \epsilon_{max}^n \sum_{k = M - M_0 - n}^{M - M_0 - 1} \left| \omega_k^{(\beta)} \right|,
\end{equation*}
where $\left| \epsilon \right| \leq \epsilon_{max}$ is used.
Thus, the proof is completed.
$\hfill\Box$

The result in Theorem \ref{th3.3} indicates that, with the new norm $\left\| \cdot \right\|_{Q,\infty}$, 
the preconditioned matrix $P_{\alpha}^{-1} \mathcal{M}_{22}$ is close to $I_{ts2}$ as $\alpha \rightarrow 0$. 
This also means that the preconditioner $P_{\alpha}$ can indeed accelerate the convergence of an iterative method.
\begin{remark}
Actually, Theorem \ref{th3.3} needs two essential properties:

(1) The quadrature weights $\left\{ \omega_k^{(\beta)} \right\}_{k = 1}^{M - M_0}$ satisfy relation \eqref{eq3.6};

(2) The spatial discretization matrix $B$ can be diagonalized.

\noindent If $B$ cannot be diagonalized, but we know that $-B$ is a nonsingular M-matrix \cite[Definition~1]{lu2015fast},
then, according to \cite[Corollary 10]{lu2015fast}, $P_\alpha$ is invertible for $\alpha \in (0,1)$ and
\begin{equation*}
\frac{\left\| P_\alpha^{-1}  - \mathcal{M}_{22}^{-1} \right\|_{\infty}}{\left\| \mathcal{M}_{22}^{-1} \right\|_{\infty}}
= \mathcal{O}(\alpha).
\end{equation*}
\end{remark}

\section{Extension to the semilinear form of Eq.~(1.1)}
\label{sec4}

In this section, we extend our method to solve the semilinear problem \cite{liao2019unconditional} of Eq.~\eqref{eq1.1}, i.e.,
\begin{equation}
\begin{cases}
\int_{0}^{t} \xi_{1 - \beta}(t - s) \partial_{s} u(x,y,s) ds = \kappa \Delta u(x,y,t) + g(u), & (x,y,t) \in \Omega \times (0,T], \\
u(x,y,t) = 0, & (x,y) \in \partial \Omega,~0 < t \leq T, \\
u(x,y,0) = u_0(x,y), & (x,y) \in \Omega,
\end{cases}
\label{eq4.0}
\end{equation}
where $g$ is a nonlinear function and nonstiff.
After discretization, we have the following nonlinear implicit scheme
\begin{equation*}
\delta_t^\beta \bm{u}^k = B \bm{u}^k + g(\bm{u}^k), \quad\mathrm{for}\quad 1 \leq k \leq M.
\end{equation*}
Then, the all-at-once system reads
\begin{equation}
\mathcal{M} \bm{u} - \bm{G}(\bm{u}) = \bm{\eta},
\label{eq4.11}
\end{equation}
where $\bm{G}(\bm{u}) = \left[ g(\bm{u}^1)^T, \cdots, g(\bm{u}^M)^T \right]^T$.
Similar to the linear case, this system is split into two subproblems:
\begin{subequations}
\begin{align}
\bm{G}_1(\tilde{\bm{u}}_1) 
&= \mathcal{M}_{11} \tilde{\bm{u}}_1 - \bm{G}(\tilde{\bm{u}}_1) - \bm{\eta}_1
= \bm{0}, \label{eq4.1a} \\
\bm{G}_2(\tilde{\bm{u}}_2) 
&= \mathcal{M}_{22} \tilde{\bm{u}}_2 - \bm{G}(\tilde{\bm{u}}_2) 
- \bm{\eta}_2 - \mathcal{M}_{21} \tilde{\bm{u}}_1 = \bm{0}. \label{eq4.1b}
\end{align}
\label{eq4.1}
\end{subequations}

In this paper, both of them are solved by a modified Newton method.
Again, to accelerate the speed of a Krylov subspace method for solving the linearized equations,
the two preconditioners $P_1$ and $P_{\alpha}$ are used.

Now, we derive our modified Newton method for solving \eqref{eq4.1}.
For the subproblem \eqref{eq4.1a}, its solution can be obtained from the following iteration process with an initial value $\tilde{\bm{u}}_1^{(0)}$
\begin{equation*}
\mathcal{M}_{11} \bm{U}_1^{\ell} = \bm{G}_1(\tilde{\bm{u}}_1^{(\ell)}), \qquad
\tilde{\bm{u}}_1^{(\ell + 1)} = \tilde{\bm{u}}_1^{(\ell)} - \bm{U}_1^{\ell}.
\end{equation*}
Then, the preconditioner $P_1$ can be directly used to accelerate solving the above equation. 
%where
%\begin{equation*}
%\mathcal{J}_1^{\ell} = \mathcal{M}_{11} - \nabla \bm{G}_1^{\ell} \quad \mathrm{with} \quad
%\nabla \bm{G}_1^{\ell} = \frac{\partial }{\partial \tilde{\bm{u}}^{1(\ell)}} \bm{G}(\tilde{\bm{u}}^{1(\ell)})
%= \mathrm{blkdiag} \left( \nabla g \left(\bm{u}^{1(\ell)} \right), \cdots, \nabla g \left( \bm{u}^{M_0 (\ell)} \right) \right).
%\end{equation*}
%Then, the preconditioner is
%\begin{equation*}
%P_{J}^{\ell} = tri(A_{11}) \otimes I_s - I_{t1} \otimes \left( B + \bar{\nabla} \bm{G}_1^{\ell} \right),
%\end{equation*}
%where
%\begin{equation*}
%\bar{\nabla} \bm{G}_1^{\ell} = \frac{1}{M_0} \sum_{n = 1}^{M_0} \nabla g \left( \bm{u}^{n(\ell)} \right).
%\end{equation*}

Similarly, the solution of \eqref{eq4.1b} can be obtained from
\begin{equation*}
\mathcal{M}_{22} \bm{U}_2^{\ell} = \bm{G}_2(\tilde{\bm{u}}_2^{(\ell)}), \qquad
\tilde{\bm{u}}_2^{(\ell + 1)} = \tilde{\bm{u}}_2^{(\ell)} - \bm{U}_2^{\ell}
\end{equation*}
with an initial value $\tilde{\bm{u}}_2^{(0)}$.
The preconditioner $P_{\alpha}$ can be used to solve this equation efficiently.
%Here
%\begin{equation*}
%\mathcal{J}_2^{\ell} = \mathcal{M}_{22} - \nabla \bm{G}_2^{\ell} \quad \mathrm{and} \quad
%\nabla \bm{G}_2^{\ell} = \frac{\partial }{\partial \tilde{\bm{u}}^{2(\ell)}} \bm{G}(\tilde{\bm{u}}^{2(\ell)})
%= \mathrm{blkdiag} \left( \nabla g \left(\bm{u}^{{M_0 + 1}(\ell)} \right), \cdots, \nabla g \left( \bm{u}^{M (\ell)} \right) \right).
%\end{equation*}
%Our preconditioner for this Jacobian equation is
%\begin{equation*}
%P_{J}^{\alpha,\ell} = A_{22}^{\alpha} \otimes I_s - I_{t2} \otimes \left( B + \bar{\nabla} \bm{G}_2^{\ell} \right),
%\end{equation*}
%where
%\begin{equation*}
%\bar{\nabla} \bm{G}_2^{\ell} = \frac{1}{M - M_0} \sum_{n = M_0 + 1}^{M} \nabla g \left( \bm{u}^{n(\ell)} \right).
%\end{equation*}
%
%Notice the diagonalization of $A_{22}^{\alpha}$, we can compute $\bm{z} = \left( P_{J}^{\alpha,\ell} \right)^{-1} \bm{v}$ in the following way:
%\begin{equation}
%\begin{cases}
%\bm{z}_1 = \left( V_\alpha^{-1} \otimes I_s \right) \bm{v}
%= \left( \mathbb{F} \otimes I_s \right) \left[ \left( \Theta_{\alpha}^{-1} \otimes I_s \right) \bm{v} \right], & \textrm{Step-(na)}, \\
%\left[ \lambda_n^{(\alpha)} I_s - \left( B + \bar{\nabla} \bm{G}_2^{\ell} \right) \right] \bm{z}_{2,n}
%= \bm{z}_{1,n}, ~1 \leq n \leq M - M_0, & \textrm{Step-(nb)}, \\
%\bm{z} = \left( V_\alpha \otimes I_s \right) \bm{z}_2
%= \left( \Theta_{\alpha} \otimes I_s \right) \left[ \left( \mathbb{F}^{*} \otimes I_s \right) \bm{z}_2 \right], & \textrm{Step-(nc)}.
%\end{cases}
%\label{eq4.2}
%\end{equation}

In the modified Newton method, the choice of the initial values $\tilde{\bm{u}}_1^{(0)}$ and $\tilde{\bm{u}}_2^{(0)}$ can affect its convergence rate.
	Thus, they should be chosen with care.
In this work, following the idea of \cite{zhao2020preconditioning},
the initial values $\tilde{\bm{u}}_1^{(0)}$ and $\tilde{\bm{u}}_2^{(0)}$ are obtained by interpolating the numerical solution of
the following linearized scheme on the coarser mesh:
\begin{equation*}
\delta_t^\beta \bm{u}^k = B \bm{u}^k + g(\bm{u}^{k - 1}), \quad\mathrm{for}\quad 1 \leq k \leq M.
\end{equation*}
Moreover, the modified Newton method terminates if $\frac{\| \bm{U}_k^{\ell} \|_2}{\| \tilde{\bm{u}}_k^{(0)} \|_2} \leq 10^{-10}~(k = 1,2)$
or the iteration number is more than $200$.
% Table 0
\begin{table}[ht]\small\tabcolsep=10.0pt
	\caption{Summary of used abbreviations.}
	\begin{center}
		\begin{tabular}{|l|l|}
			\hline
			Symbol & Explanation \\
			\hline
			BFSM & The BFS method is used to solve \eqref{eq2.4} or \eqref{eq4.1}. \\
			$\mathcal{I}$ & The BiCGSTAB method is used when solving \eqref{eq2.4} or \eqref{eq4.1}. \\
			$\mathcal{P}$ & The PBiCGSTAB method is used when solving \eqref{eq2.4} or \eqref{eq4.1}. \\
			$\mathrm{Iter}$ & $\mathrm{Iter} = \left( \mathrm{Iter}(1), \mathrm{Iter}(2) \right)$, 
			where $\mathrm{Iter}(1)$ ($\mathrm{Iter}(2)$) is the number of iterations required \\
			& for solving Eq.~\eqref{eq3.2a} (Eq.~\eqref{eq3.2b}). \\
			$\mathrm{Iter1}$ & $\mathrm{Iter1} = \frac{1}{M} \sum\limits_{k = 1}^{M} \mathrm{Iter1}(k)$, 
			where $\mathrm{Iter1}(k)$ is the number of iterations of the modified \\
			& Newton method in the $k$th step of BFSM for solving \eqref{eq4.11}. \\
			$\mathrm{Iter}_{O}$ & $\mathrm{Iter}_{O} = \left( \mathrm{Iter}_{O}(1), \mathrm{Iter}_{O}(2) \right)$, 
			where $\mathrm{Iter}_{O}(1)$ ($\mathrm{Iter}_{O}(2)$) is the number of iterations of \\
			& the modified Newton method for solving Eq.~\eqref{eq4.1a} (Eq.~\eqref{eq4.1b}). \\
			%$\mathrm{Iter}_{O2}$ & The number of iterations of the modified Newton's method for solving  \\
			$\mathrm{Iter}_{I}$ & $\mathrm{Iter}_{I} = \left( \mathrm{Iter}_{I}(1), \mathrm{Iter}_{I}(2) \right)$,
			where $\mathrm{Iter}_{I}(1) = \frac{1}{\mathrm{Iter}_{O}(1)} \sum\limits_{k = 1}^{\mathrm{Iter}_{O}(1)} \mathrm{Iter}_{I1}(k)$ and \\
			& $\mathrm{Iter}_{I}(2) = \frac{1}{\mathrm{Iter}_{O}(2)} \sum\limits_{k = 1}^{\mathrm{Iter}_{O}(2)} \mathrm{Iter}_{I2}(k)$. \\
			$\mathrm{Iter}_{I1}(k)$ & The number of iterations required by the (P)BiCGSTAB method in the $k$th step \\
			& of the modified Newton method for solving Eq.~\eqref{eq4.1a}. \\
			$\mathrm{Iter}_{I2}(k)$ & The number of iterations required by the (P)BiCGSTAB method in the $k$th step \\
			& of the modified Newton method for solving Eq.~\eqref{eq4.1b}. \\
			%$\mathrm{Iter}_{I2}$ & $\mathrm{Iter}_{I2} = \frac{1}{\mathrm{Iter}_{O2}} \sum\limits_{k = 1}^{\mathrm{Iter}_{O2}}
			%$\mathrm{Iter}_{I2}(k)$, where the meaning of $\mathrm{Iter}_{I2}(k)$ is similar to $\mathrm{Iter}_{I1}(k)$ \\
			$\textrm{Time}$ & Total CPU time in seconds. \\
			-- & The data is not obtained in 8 hours. \\
			\dag & The BiCGSTAB method needs more than 1000 iterations to reach the desired tolerance. \\
			\hline
		\end{tabular}
		\label{tab0}
	\end{center}
\end{table}
\section{Numerical experiments}
\label{sec5}
In this section, three examples are reported to show the performance
of our strategy proposed in Sections \ref{sec3} and \ref{sec4}.
The (P)BiCGSTAB method for solving \eqref{eq3.2} terminates
if the relative residual error satisfies
$\frac{\| \bm{r}^{(k)} \|_2}{\| \bm{r}^{(0)} \|_2} \leq 10^{-9}$ (for the nonlinear case we chose $10^{-6}$)
or the iteration number is more than $1000$, where $\bm{r}^{(k)}$ denotes residual vector in the $k$th iteration.
The initial guess is chosen as the zero vector. 
In our experiments, we set $N = N_x = N_y$, $T_0 = 2^{-r}$, 
$M_0 = \left\lceil \frac{r}{2^r - 1 + r} M \right\rceil$ and $\alpha = \min \{10^{-4}, 0.5 \tilde{\tau} \}$.
%It is easy to check that $\tilde{\tau} \geq \tau_{M_0}$. 
All of the symbols shown in Table \ref{tab0} will appear later.

All experiments are carried out via MATLAB 2018b on a Windows 10 (64 bit) PC with the configuration:
Intel(R) Core(TM) i7-8700k CPU 3.20 GHz and 16 GB RAM.
% Table 1
\begin{table}[ht]\footnotesize\tabcolsep=8.0pt
	\begin{center}
		\caption{Results of various methods for $M = N$ for Example 1.}
		\centering
		\begin{tabular}{ccccccc}
			\hline
			%& & \multicolumn{5}{c}{$r = 2$} & \multicolumn{5}{c}{$r = 3$} \\
			%[-2pt] \cmidrule(lr){3-7} \cmidrule(lr){8-12} \\ [-11pt]
			& & \rm{BFSM} & \multicolumn{2}{c}{$\mathcal{I}$} & \multicolumn{2}{c}{$\mathcal{P}$} \\
			[-2pt] \cmidrule(lr){4-5} \cmidrule(lr){6-7} \\ [-11pt]
			$(\beta, r)$ & $N$ & $\mathrm{Time}$ & $\mathrm{Iter}$ & $\mathrm{Time}$ & $\mathrm{Iter}$ & $\mathrm{Time}$ \\
			\hline
			(0.1, 2) & 32 & 0.037 & (37.0, 41.0) & 0.251 & (3.0, 1.0) & 0.032 \\ 
			& 64 & 0.328 & (67.0, 71.0) & 1.645 & (4.0, 1.0) & 0.182 \\ 
			& 128 & 4.190 & (136.0, 145.0) & 36.962 & (4.0, 1.0) & 1.399 \\
			& 256 & 44.256 & (265.0, 303.0) & 906.017 & (5.0, 1.0) & 18.111 \\
			& 512 & 1108.741 & -- & -- & (5.0, 1.0) & 224.922 \\
			\\
			(0.5, 2) & 32 & 0.039 & (27.0, 31.0) & 0.190 & (3.0, 1.0) & 0.033 \\
			& 64 & 0.326 & (55.0, 77.0) & 1.706 & (5.0, 2.0) & 0.263 \\ 
			& 128 & 3.937 & (128.0, 184.0) & 44.789 & (7.0, 2.0) & 2.425 \\
			& 256 & 44.382 & (308.0, 401.0) & 1178.911 & (10.0, 2.0) & 34.433 \\ 
			& 512 & 1115.207 & -- & -- & (14.0, 2.0) & 551.918 \\ 
			\\
			(0.9, 2) & 32 & 0.050 & (32.0, 24.0) & 0.154 & (3.0, 2.0) & 0.042 \\ 
			& 64 & 0.344 & (64.0, 71.0) & 1.632 & (4.0, 2.0) & 0.230 \\
			& 128 & 3.890 & (177.0, 213.0) & 53.116 & (5.0, 2.0) & 2.050 \\
			& 256 & 44.434 & (579.0, 638.0) & 1915.291 & (6.0, 2.0) & 25.939 \\
			& 512 & 1111.531 & -- & -- & (9.0, 2.0) & 434.259 \\
			\\        
			(0.1, 3) & 32 & 0.038 & (36.0, 44.0) & 0.338 & (3.0, 1.0) & 0.032 \\
			& 64 & 0.323 & (63.0, 73.0) & 2.617 & (3.0, 1.0) & 0.162 \\
			& 128 & 3.894 & (139.0, 148.0) & 51.725 & (4.0, 1.0) & 1.392 \\
			& 256 & 44.510 & (291.0, 312.0) & 1156.569 & (4.0, 1.0) & 15.548 \\
			& 512 & 1115.510 & -- & -- & (5.0, 1.0) & 217.845 \\
			\\
			(0.5, 3) & 32 & 0.036 & (33.0, 36.0) & 0.277 & (3.0, 1.0) & 0.032 \\
			& 64 & 0.320 & (56.0, 86.0) & 3.013 & (4.0, 2.0) & 0.255 \\
			& 128 & 3.911 & (127.0, 197.0) & 66.950 & (5.0, 2.0) & 2.194 \\
			& 256 & 44.728 & (267.0, 422.0) & 1510.010 & (7.0, 2.0) & 27.578 \\
			& 512 & 1108.584 & -- & -- & (10.0, 2.0) & 442.688 \\
			\\
			(0.9, 3) & 32 & 0.035 & (48.0, 28.0) & 0.224 & (2.0, 2.0) & 0.038 \\
			& 64 & 0.323 & (121.0, 85.0) & 3.102 & (3.0, 2.0) & 0.225 \\
			& 128 & 3.882 & (386.0, 256.0) & 93.875 & (4.0, 2.0) & 2.042 \\
			& 256 & 45.184 & \dag & \dag & (5.0, 2.0) & 24.569 \\
			& 512 & 1111.146 & -- & -- & (6.0, 2.0) & 373.949 \\
			\hline
		\end{tabular}
		\label{tab1}
	\end{center}
\end{table}

\noindent \textbf{Example 1.}
%Considering Eq. (\ref{eq1.1}) with $\Omega = [0,\pi]$, $T = 1$
%and the source term
%\begin{equation*}
%f(x,t) = \left[ \xi_{1 + \sigma - \beta}(t) - \kappa \xi_{1 + \sigma}(t) \right] \sin x,
%\end{equation*}
%where $\kappa = 2$ and $\sigma = 2.2 - \beta$.
%The exact solution is $u(x,t) = \xi_{1 + \sigma}(t) \sin x$.
In this example, the subdiffusion problem \eqref{eq1.1} is considered on $\Omega = [-4,10] \times [-4,10]$ with
$T = 1$ and the source term
\begin{equation*}
\begin{split}
f(x,y,t) = & \frac{\xi_{1 + \sigma - \beta}(t)}{\sqrt{2 \pi}} 
\left[\exp \left(-\frac{x^2 + y^2}{2} \right) + \exp \left(-\frac{(x - 3)^2 + (y - 3)^2}{2} \right) \right]
- \kappa \frac{1 + \xi_{1 + \sigma}(t)}{\sqrt{2 \pi}} \times \\
& \left\{ \left( x^2 + y^2 - 2 \right) \exp \left(-\frac{x^2 + y^2}{2} \right) 
+ \left[ (x - 3)^2 + (y - 3)^2 - 2 \right] \exp \left(-\frac{(x - 3)^2 + (y - 3)^2}{2} \right) \right\},
\end{split}
\end{equation*}
where $\kappa = 1$ and $\sigma = 2.2 - \beta$.
For the above choice, the exact solution is 
\begin{equation*}
u(x,y,t) = \frac{1 + \xi_{1 + \sigma}(t)}{\sqrt{2 \pi}} \left[\exp \left(-\frac{x^2 + y^2}{2} \right) + \exp \left(-\frac{(x - 3)^2 + (y - 3)^2}{2} \right) \right].
\end{equation*}
% Figure 2
\begin{figure}[ht]
	\centering
	\subfigure[Eigenvalues of $\mathcal{M}_{11}$]
	{\includegraphics[width=2.2in,height=2.2in]{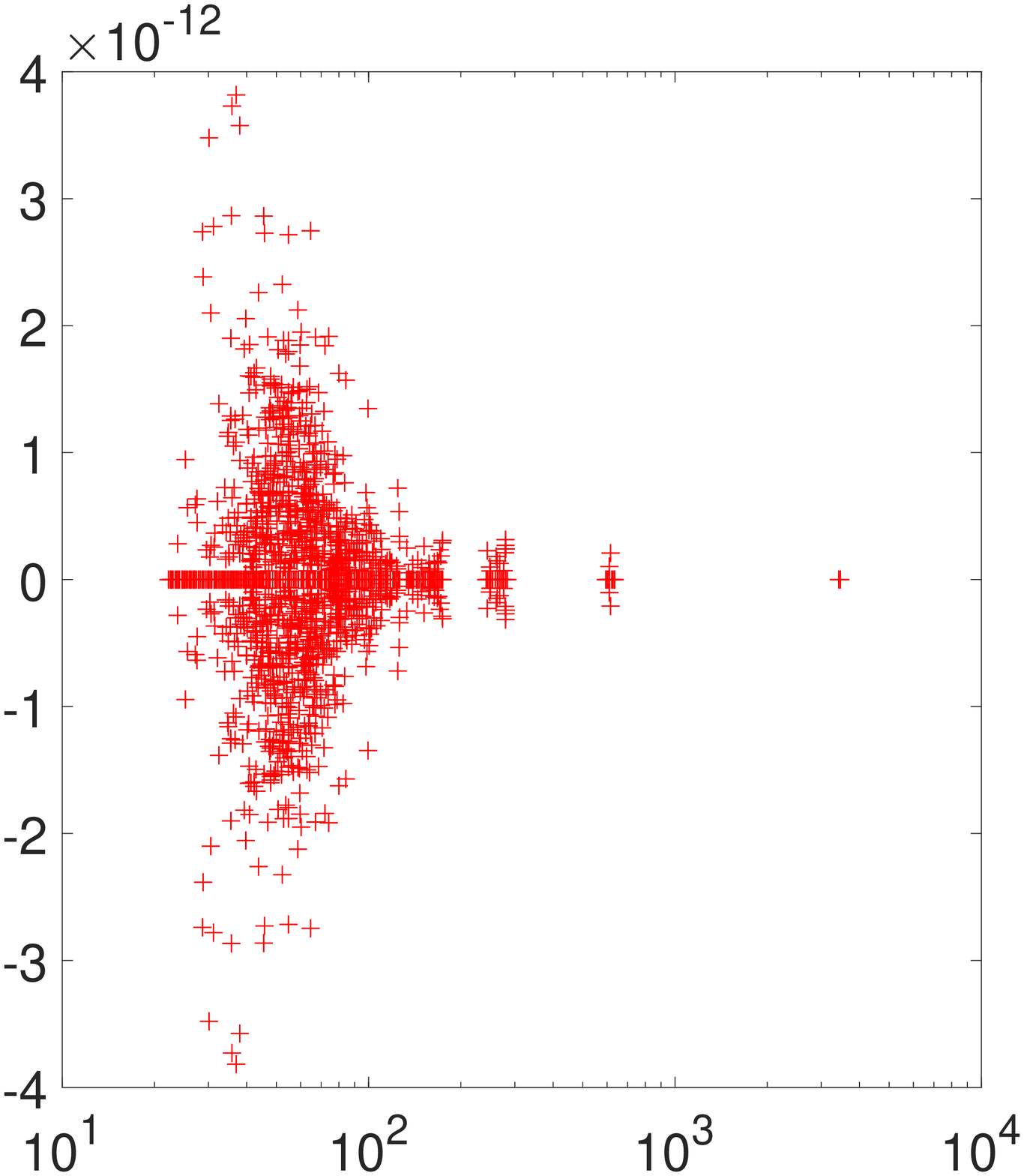}}\hspace{10mm}
	\subfigure[Eigenvalues of $P_{1}^{-1} \mathcal{M}_{11}$]
	{\includegraphics[width=2.2in,height=2.2in]{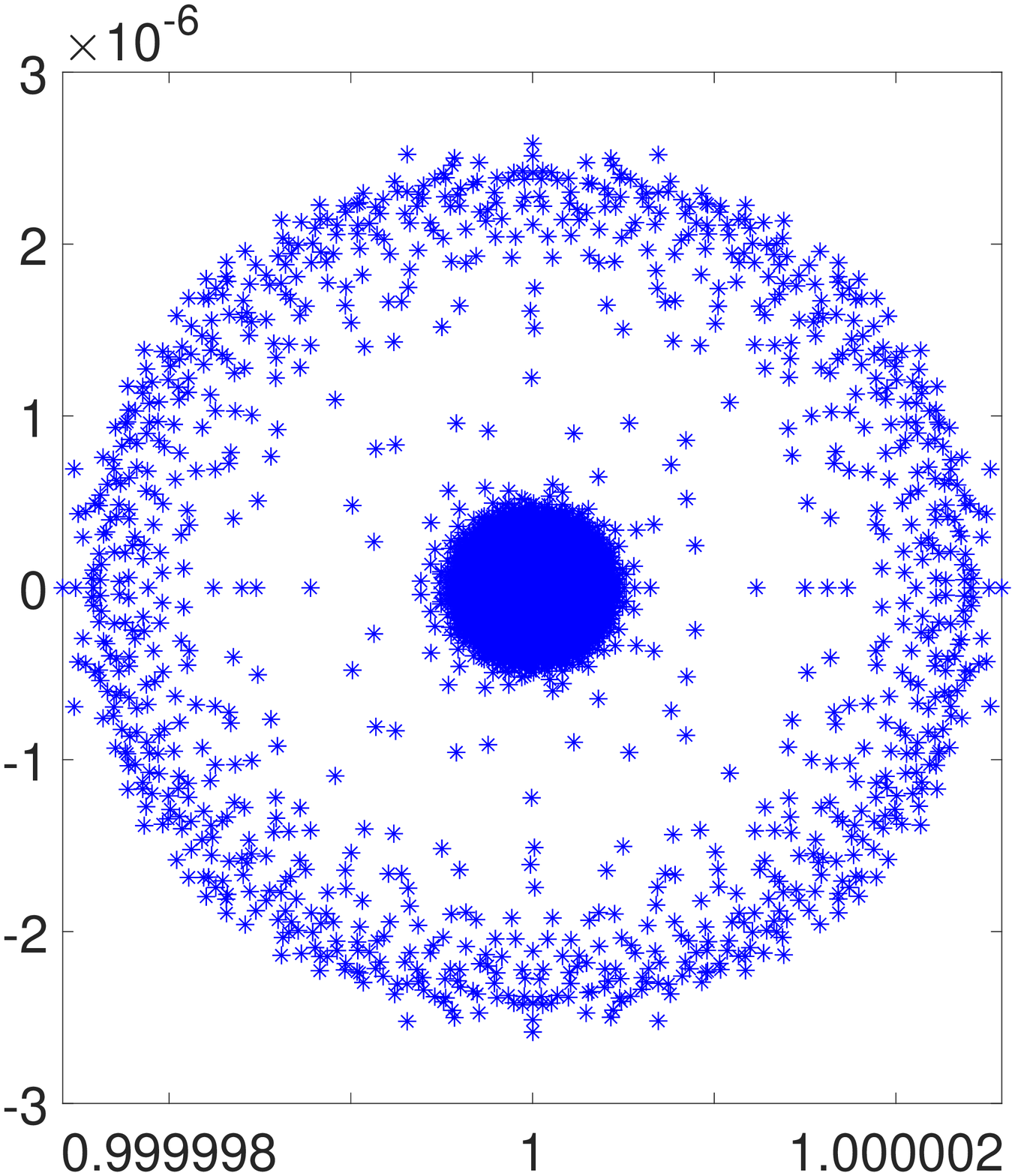}}\\ \hspace{5mm}
	\subfigure[Eigenvalues of $\mathcal{M}_{22}$]
	{\includegraphics[width=2.2in,height=2.2in]{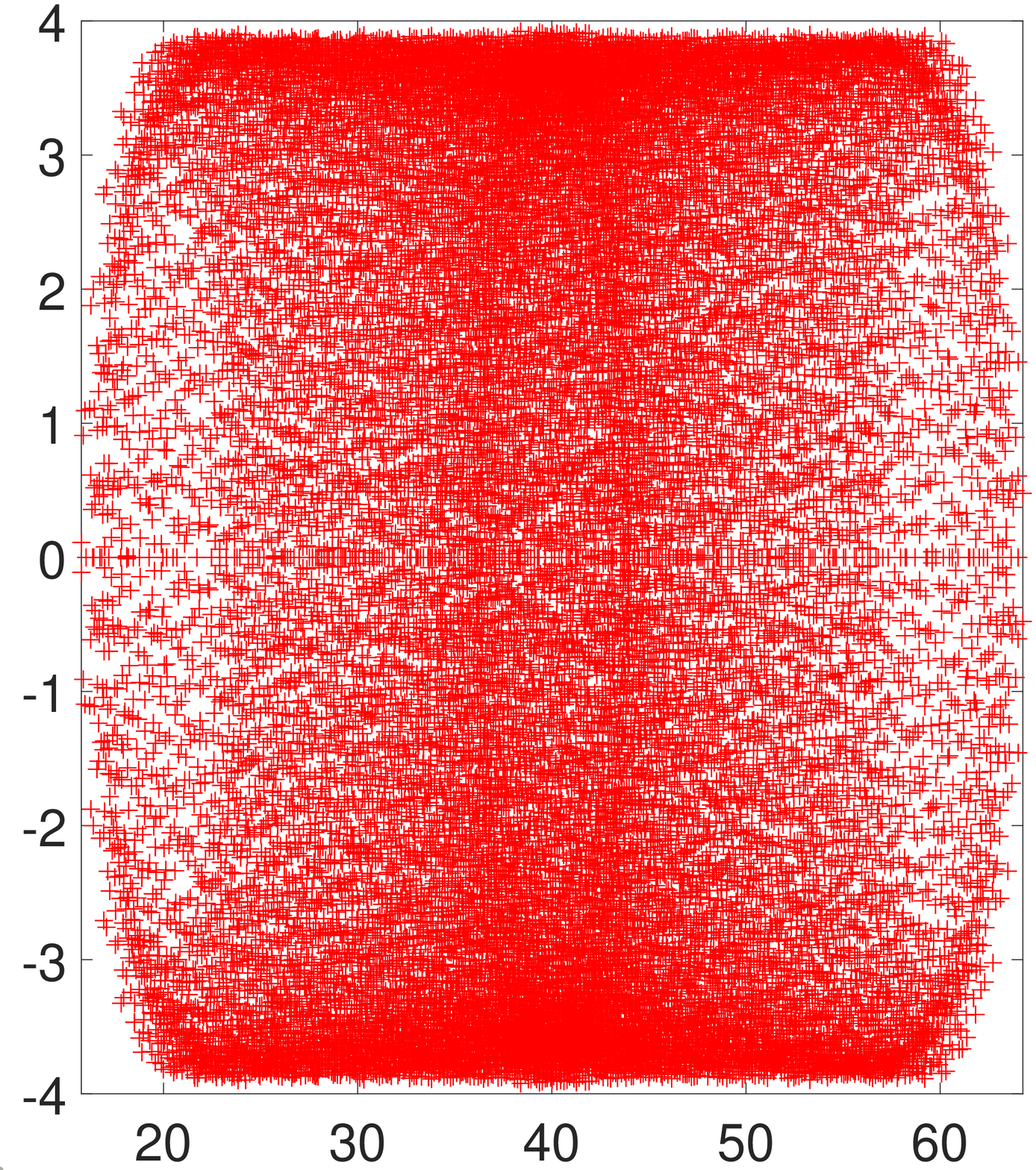}}\hspace{7mm}
	\subfigure[Eigenvalues of $P_{\alpha}^{-1} \mathcal{M}_{22}$]
	{\includegraphics[width=2.5in,height=2.3in]{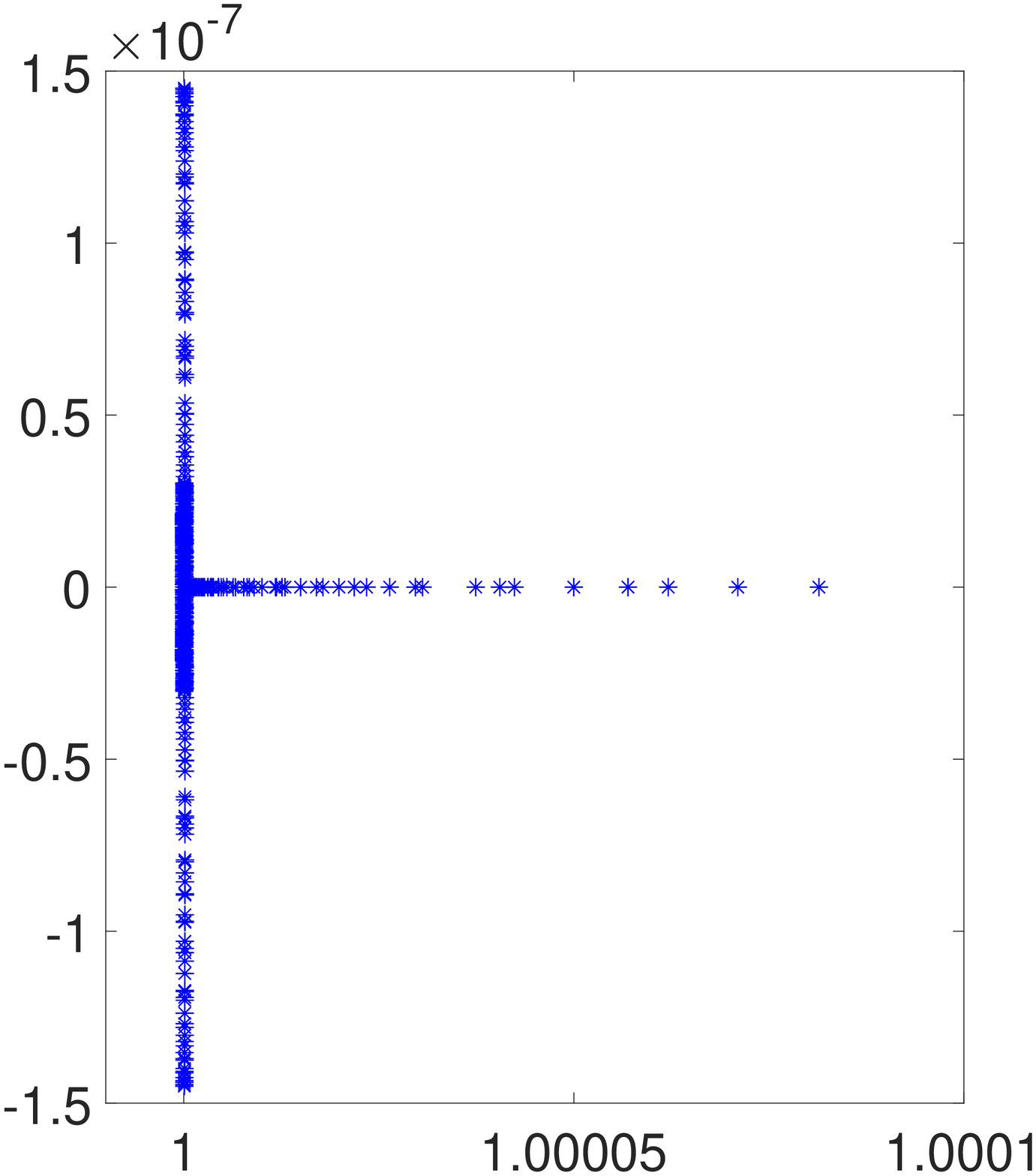}}
	\caption{Spectra of $\mathcal{M}_{11}$, $P_{1}^{-1} \mathcal{M}_{11}$, $\mathcal{M}_{22}$ and $P_{\alpha}^{-1} \mathcal{M}_{22}$, 
		for $(\beta, r) = (0.9, 3)$ and $M = N = 32$ in Example 1.}
	\label{fig2}
\end{figure}

In Table \ref{tab1}, the CPU time of method $\mathcal{P}$ is the smallest one among the three tested methods.
Comparing the number $\mathrm{Iter}$ of methods $\mathcal{I}$ and $\mathcal{P}$,
it can be found that our preconditioners $P_1$ and $P_{\alpha}$ are very efficient.
We also notice that the number $\mathrm{Iter}$ of our method (i.e., $\mathcal{P}$) is not strongly influenced by the mesh size.
Moreover, Fig.~\ref{fig2} shows the spectra of $\mathcal{M}_{11}$, $P_{1}^{-1} \mathcal{M}_{11}$, $\mathcal{M}_{22}$ and $P_{\alpha}^{-1} \mathcal{M}_{22}$ with $(\beta, r) = (0.9, 3)$ and $M = N = 32$.
It should be mentioned that the effect of the clustered eigenvalues on the convergence of the (P)BiCGSTAB method can be not so crucial \cite{Greenbaum1997iterative}.
% Table 2
\begin{table}[ht]\footnotesize\tabcolsep=6.0pt
	\begin{center}
		\caption{Results of various methods for $M = N$ for Example 2.}
		\centering
		\begin{tabular}{cccccccccc}
			\hline
			%& & \multicolumn{5}{c}{$r = 2$} & \multicolumn{5}{c}{$r = 3$} \\
			%[-2pt] \cmidrule(lr){3-7} \cmidrule(lr){8-12} \\ [-11pt]
			& & \multicolumn{2}{c}{\rm{BFSM}} & \multicolumn{3}{c}{$\mathcal{I}$} & \multicolumn{3}{c}{$\mathcal{P}$} \\
			[-2pt] \cmidrule(lr){3-4} \cmidrule(lr){5-7} \cmidrule(lr){8-10} \\ [-11pt]
			$(\beta, r)$ & $N$ & $\mathrm{Iter1}$ & $\mathrm{Time}$ & $\mathrm{Iter}_{O}$ & $\mathrm{Iter}_{I}$ & $\mathrm{Time}$ 
			& $\mathrm{Iter}_{O}$ & $\mathrm{Iter}_{I}$ & $\mathrm{Time}$ \\
			\hline
			 (0.1, 2) & 32 & 9.8 & 0.263 & (10.0, 10.0) & (66.9, 55.7) & 3.409 & (10.0, 10.0) & (2.0, 1.0) & 0.320 \\
			          & 64 & 9.0 & 2.193 & (10.0, 10.0) & (135.9, 100.9) & 24.137 & (10.0, 10.0) & (2.1, 1.0) & 1.283 \\
			          & 128 & 9.0 & 25.604 & (10.0, 10.0) & (271.2, 219.3) & 579.351 & (10.0, 10.0) & (2.8, 1.0) & 10.991 \\
			          & 256 & 8.1 & 238.001 & (10.0, 10.0) & (607.2, 457.3) & 14827.553 & (10.0, 10.0) & (2.9, 1.0) & 138.210 \\
			          & 512 & 8.0 & 3115.653 & -- & -- & -- & (10.0, 10.0) & (3.0, 1.0) & 1789.639 \\
			 \\
			 (0.5, 2) & 32 & 7.7 & 0.207 & (8.0, 10.0) & (61.0, 65.2) & 3.869 & (8.0, 10.0) & (2.4, 1.0) & 0.303 \\
			          & 64 & 7.1 & 1.772 & (8.0, 10.0) & (139.9, 135.4) & 30.318 & (8.0, 10.0) & (3.5, 1.0) & 1.468 \\
			          & 128 & 6.7 & 19.401 & (8.0, 10.0) & (318.4, 295.3) & 742.383 & (8.0, 10.0) & (4.5, 1.0) & 12.320 \\
			          & 256 & 6.0 & 181.785 & (8.0, 10.0) & (701.5, 642.7) & 19232.301 & (8.0, 10.0) & (6.1, 1.0) & 181.685 \\
			          & 512 & 5.7 & 2469.990 & -- & -- & -- & (8.0, 10.0) & (8.0, 1.0) & 2628.717 \\
			 \\
			 (0.9, 2) & 32 & 6.4 & 0.175 & (6.0, 9.0) & (42.5, 57.7) & 2.993 & (6.0, 9.0) & (2.0, 1.0) & 0.262 \\
			          & 64 & 5.6 & 1.427 & (6.0, 9.0) & (104.2, 152.3) & 28.987 & (6.0, 9.0) & (2.7, 1.0) & 1.107 \\
			          & 128 & 4.9 & 14.642 & (6.0, 9.0) & (279.5, 418.4) & 873.289 & (6.0, 9.0) & (3.2, 1.0) & 9.044 \\
			          & 256 & 4.7 & 146.312 & (6.0, 9.0) & (841.2, 856.4) & 22173.364 & (6.0, 9.0) & (4.0, 1.0) & 123.540 \\
			          & 512 & 4.0 & 1978.387 & -- & -- & -- & (6.0, 9.0) & (6.2, 1.0) & 1842.001 \\
			 \\
			 (0.1, 3) & 32 & 9.8 & 0.269 & (9.0, 10.0) & (66.9, 59.8) & 4.424 & (9.0, 10.0) & (2.0, 1.0) & 0.253 \\
			          & 64 & 9.0 & 2.249 & (9.0, 10.0) & (135.4, 110.1) & 39.107 & (9.0, 10.0) & (2.0, 1.0) & 1.388 \\
			          & 128 & 8.9 & 25.403 & (9.0, 10.0) & (276.1, 231.6) & 817.984 & (9.0, 10.0) & (2.0, 1.0) & 11.298 \\
			          & 256 & 8.2 & 237.566 & (9.0, 10.0) & (595.4, 488.2) & 18366.031 & (9.0, 10.0) & (2.6, 1.0) & 133.166 \\
			          & 512 & 8.0 & 3031.510 & -- & -- & -- & (9.0, 10.0) & (2.9, 1.0) & 1795.297 \\
			 \\
			 (0.5, 3) & 32 & 7.6 & 0.215 & (7.0, 10.0) & (51.6, 66.5) & 4.882 & (7.0, 10.0) & (2.0, 1.0) & 0.228 \\
			          & 64 & 6.9 & 1.759 & (7.0, 10.0) & (116.4, 140.8) & 47.844 & (7.0, 10.0) & (2.9, 1.0) & 1.463 \\
			          & 128 & 6.7 & 19.198 & (7.0, 10.0) & (258.3, 296.9) & 984.842 & (7.0, 10.0) & (3.6, 1.0) & 12.189 \\
			          & 256 & 5.9 & 174.579 & (7.0, 10.0) & (601.3, 648.1) & 23378.635 & (7.0, 10.0) & (5.1, 1.0) & 152.465 \\
			          & 512 & 5.7 & 2410.630 & -- & -- & -- & (7.0, 10.0) & (6.4, 1.0) & 2161.753 \\
			 \\
			 (0.9, 3) & 32 & 6.3 & 0.180 & (6.0, 9.0) & (58.5, 60.8) & 3.999 & (6.0, 9.0) & (2.0, 1.0) & 0.199 \\ 
			          & 64 & 5.7 & 1.450 & (6.0, 9.0) & (144.8, 153.8) & 46.973 & (6.0, 9.0) & (2.0, 1.0) & 1.171 \\
			          & 128 & 4.8 & 14.497 & (6.0, 9.0) & (431.0, 430.8) & 1294.543 & (6.0, 9.0) & (2.3, 1.0) & 9.938 \\
			          & 256 & 4.7 & 143.793 & \dag & \dag & \dag & (6.0, 9.0) & (3.3, 1.0) & 119.759 \\
			          & 512 & 3.9 & 1913.440 & -- & -- & -- & (6.0, 9.0) & (3.8, 1.0) & 1626.390 \\
			\hline
		\end{tabular}
		\label{tab2}
	\end{center}
\end{table}

\noindent \textbf{Example 2.}
We consider the two-dimensional time fractional Fisher equation. More precisely,
in Eq.~\eqref{eq4.0}, we choose $\Omega = [0,\pi] \times [0,\pi]$, $T = 1$,
the diffusion coefficient $\kappa = 1$, the nonlinear term $g(u) = u (1 - u)$ 
and the initial value $u_0(x,y) = \sin x \sin y$.

\begin{figure}[t]
	\centering
	\subfigure[Eigenvalues of $\mathcal{M}_{11}$]
	{\includegraphics[width=2.2in,height=2.2in]{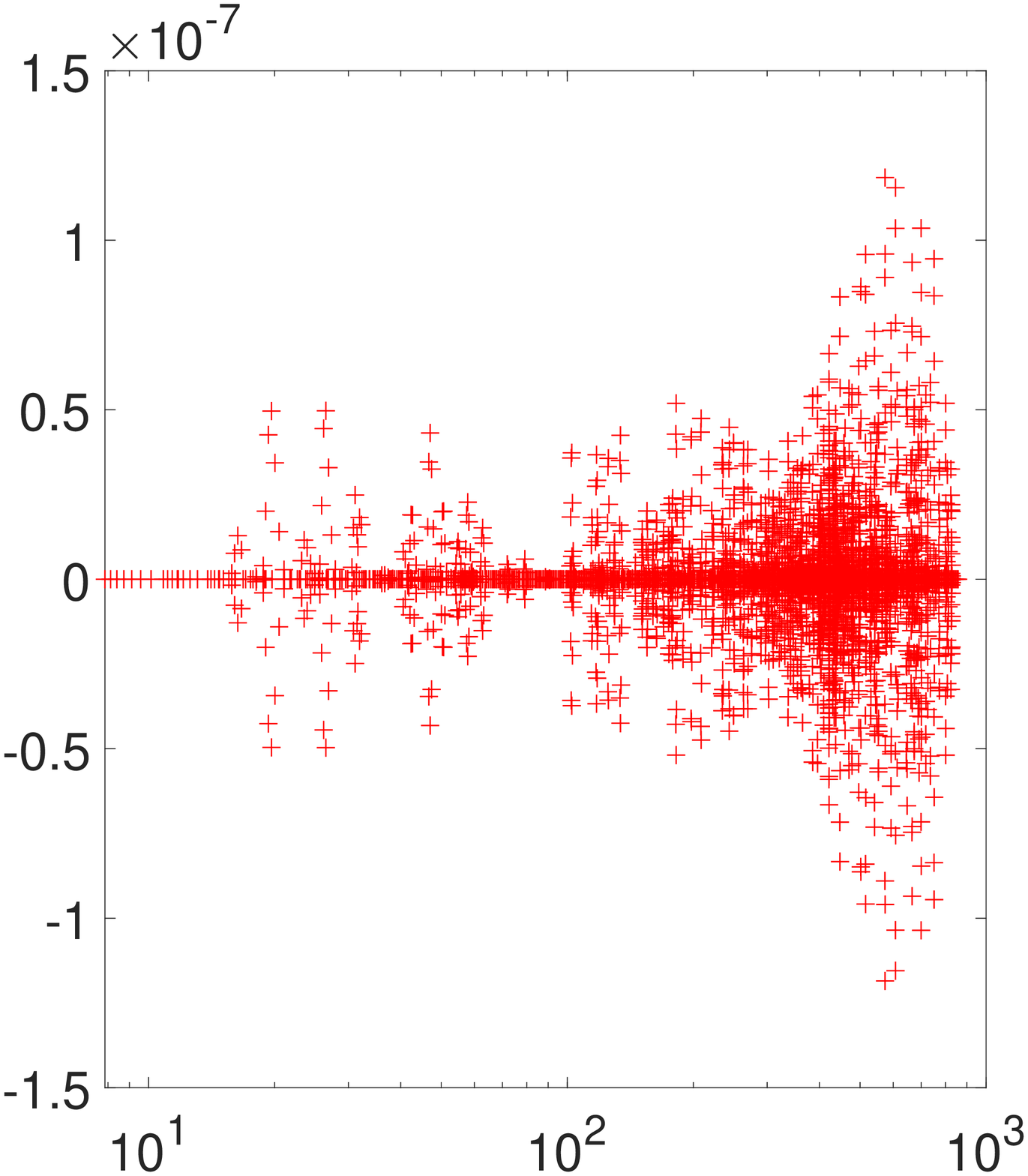}}\hspace{10mm}
	\subfigure[Eigenvalues of $P_{1}^{-1} \mathcal{M}_{11}$]
	{\includegraphics[width=2.2in,height=2.2in]{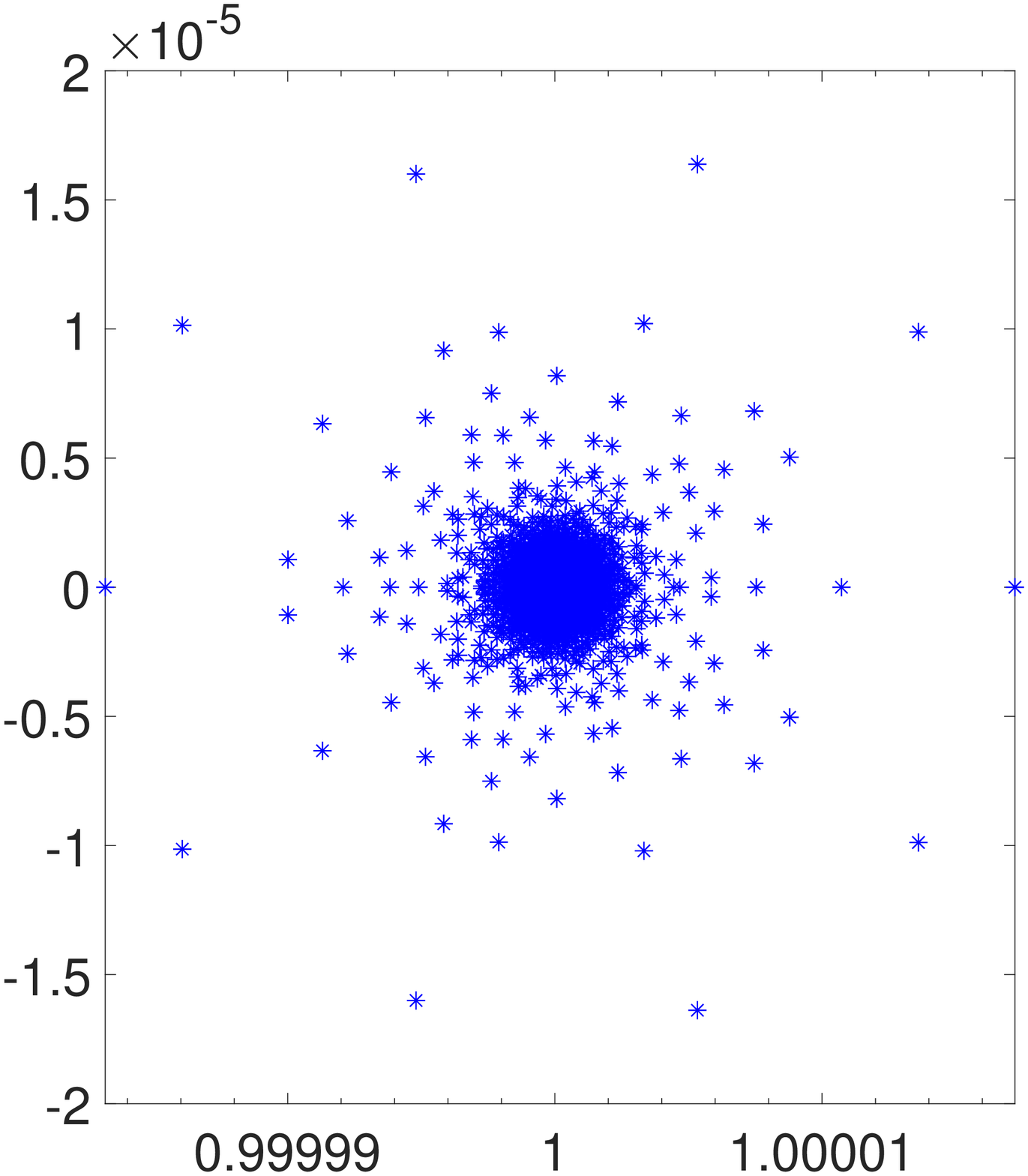}}\\ \hspace{4mm}
	\subfigure[Eigenvalues of $\mathcal{M}_{22}$]
	{\includegraphics[width=2.3in,height=2.2in]{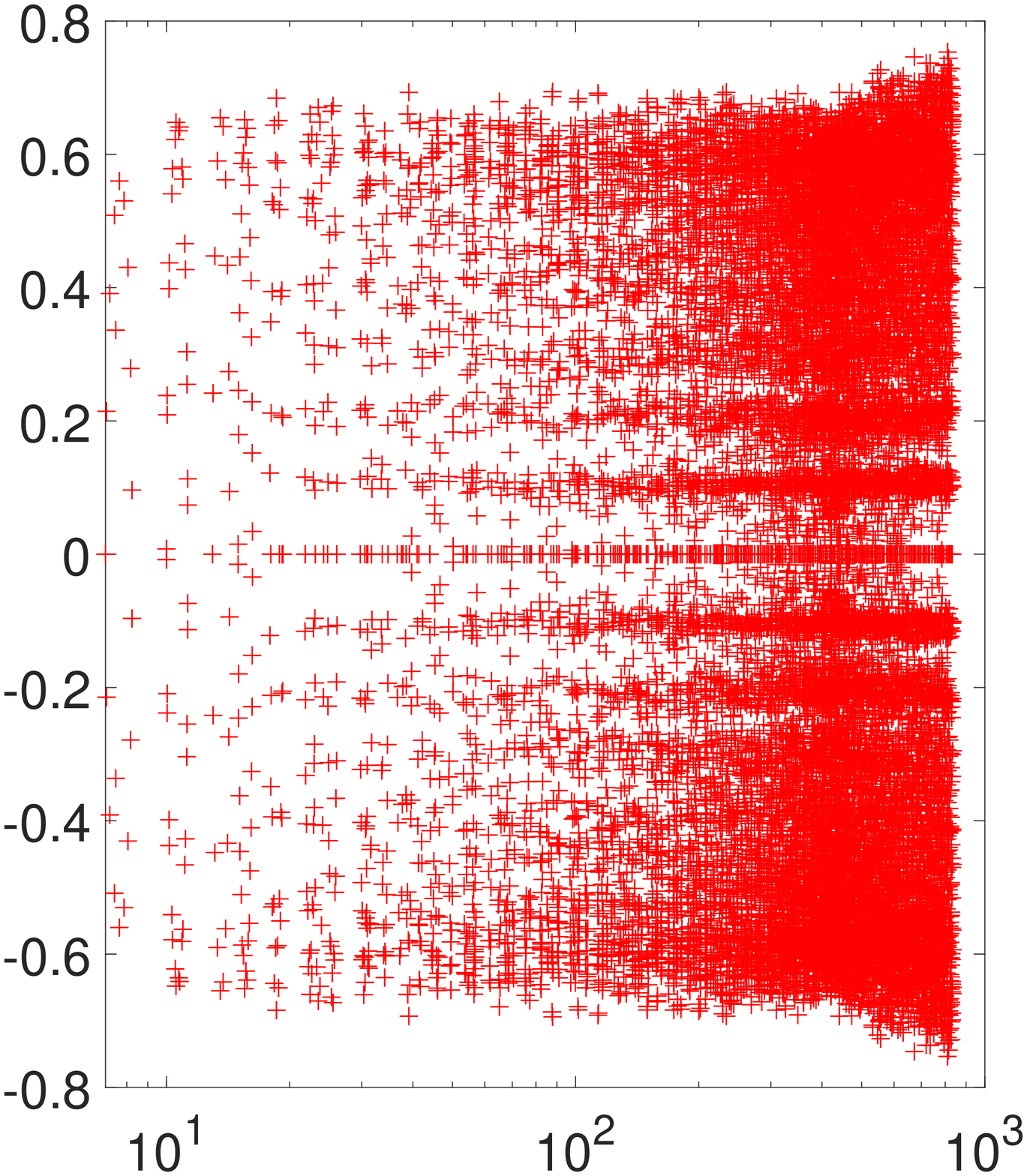}}\hspace{10mm}
	\subfigure[Eigenvalues of $P_{\alpha}^{-1} \mathcal{M}_{22}$]
	{\includegraphics[width=2.3in,height=2.3in]{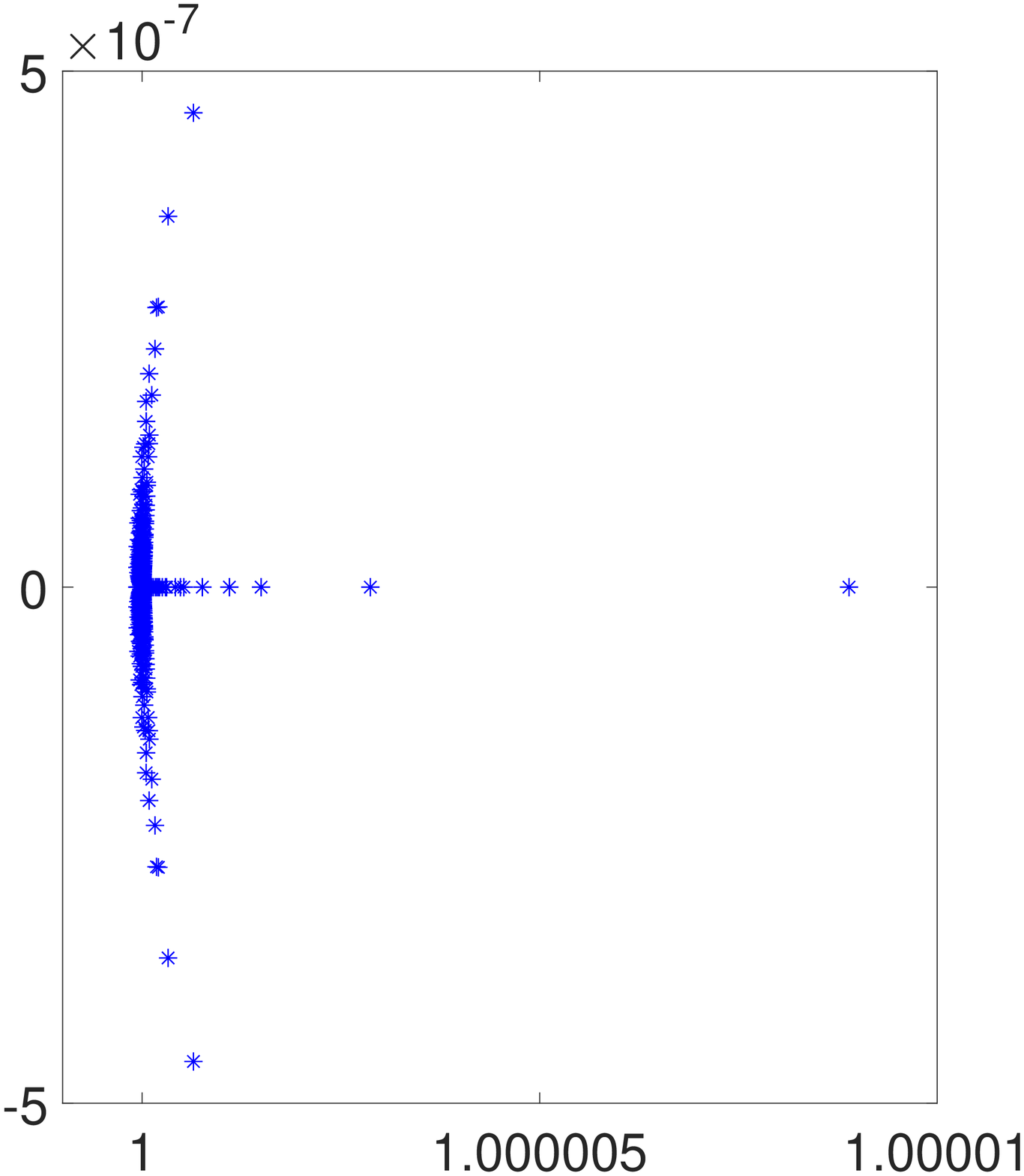}}
	\caption{Spectra of $\mathcal{M}_{11}$, $P_{1}^{-1} \mathcal{M}_{11}$, $\mathcal{M}_{22}$ and $P_{\alpha}^{-1} \mathcal{M}_{22}$, 
		for $(\beta, r) = (0.5, 2)$ and $M = N = 32$ in Example 2.}
	\label{fig3}
\end{figure}

In Table \ref{tab2}, the CPU time and the numbers of iterations of the methods BFSM, $\mathcal{I}$ and $\mathcal{P}$ for solving the nonlinear problem are reported.
Compared with the BFSM method, our method indeed reduces the CPU time except for some cases.
For these unsatisfied cases, although the CPU times required by our method $\mathcal{P}$ are larger than the BFSM method, 
our method still has a potential advantage in terms of parallel computing. 
Fig.~\ref{fig3} shows the spectra of $\mathcal{M}_{11}$, $P_{1}^{-1} \mathcal{M}_{11}$, $\mathcal{M}_{22}$ and $P_{\alpha}^{-1} \mathcal{M}_{22}$
with $(\beta, r) = (0.5, 2)$ and $M = N = 32$.
It is clearly seen from Fig.~\ref{fig3} that all eigenvalues of the preconditioned matrices 
$P_{1}^{-1} \mathcal{M}_{11}$ and $P_{\alpha}^{-1} \mathcal{M}_{22}$ are clustered around $1$.

% Figure 3

\section{Concluding remarks}
\label{sec6}

A parallel preconditioning technique is proposed to solve the all-at-once system \eqref{eq2.4} with variable time steps arising from subdiffusion equations \eqref{eq1.1}.
Firstly, we split the time interval $[0,T]$ into two parts $[0,T_0]$ and $[T_0,T]$.
Then, we use the graded $L1$ scheme to approximate \eqref{eq1.1} in $[0,T_0]$, while the uniform one is applied in $[T_0,T]$.
Secondly, our all-at-once system \eqref{eq2.4} is derived based on this decomposition.
Thanks to the local Toeplitz structure of the time discretization matrix $A$,
the solution of Eq.~\eqref{eq2.4} can be obtained by solving \eqref{eq3.2}.
Two preconditioners $P_1$ and $P_{\alpha}$ are proposed to accelerate obtaining the solution of Eq.~\eqref{eq3.2}.
Some properties of these two preconditioners are also analyzed.
In Section \ref{sec4}, we extend our technique the nonlinear subdiffusion problem \eqref{eq4.0}.
Finally, numerical experiments are reported that show the performance of our preconditioning technique.
It is worth mentioning that the CPU time required by the method $\mathcal{P}$ can be further reduced since it is suitable for parallel computing.

In this work, we consider the nonlinear function $g$ to be nonstiff. 
If $g$ is stiff, we suggest using Newton's method to solve \eqref{eq4.1}.
For this case, our preconditioners need some modifications as proposed in \cite[Section 3]{gu2020parallel} to make them more efficient.
Another benefit of such modifications is that the new preconditioners are still suitable for parallel computing.
In our future work, we will study the all-at-once system with a space discretization matrix $B$ being indefinite.

\section*{Acknowledgments}
\addcontentsline{toc}{section}{Acknowledgments}
\label{sec7}

\textit{This research is supported by the National Natural Science Foundation
of China (No.~11801463) and the Applied Basic Research Program of Sichuan Province (No.~2020YJ0007).
The first author is also supported by the China Scholarship Council.
We would like to express our sincere thanks to the referees for insightful comments and invaluable suggestions
that greatly improved the representation of this paper.
}
\section*{Declarations}
\textbf{Conflict of interest} The authors declare that they have no competing interests.

%\section*{References}
%\addcontentsline{toc}{section}{References}
\bibliography{Ref}

\end{document}